\documentclass[12pt]{amsart}
\usepackage{amsmath,amsthm,amssymb,amsfonts,amscd,color}
\usepackage{epsfig}
\usepackage{xypic}
\usepackage[colorlinks, citecolor=blue]{hyperref}

\numberwithin{equation}{section}

\theoremstyle{plain}
\newtheorem{theorem}{Theorem}

\newtheorem{corollary}[theorem]{Corollary}
\newtheorem{proposition}[theorem]{Proposition}
\newtheorem*{theorem*}{Theorem}
\newtheorem*{conjecture*}{Conjecture}

\theoremstyle{definition}
\newtheorem{remark}[theorem]{Remark}
\newtheorem{example}[theorem]{Example}
\newtheorem*{definition}{Definition}

\newcommand{\CC}{{\mathbb C}}

\newcommand{\PP}{{\mathbb P}}

\newcommand{\RR}{{\mathbb R}}
\newcommand{\ZZ}{{\mathbb Z}}

\newcommand{\DD}{{\mathbb D}}

\newcommand{\TT}{{\mathbb T}}

\newcommand{\calO}{{\mathcal O}}

\newcommand{\calC}{{\mathcal C}}

\newcommand{\calE}{{\mathcal E}}

\newcommand{\calR}{{\mathcal R}}
\newcommand{\calV}{{\mathcal V}}
\newcommand{\eps}{\varepsilon}

\newcommand{\xx}{\bf x}
\newcommand{\zz}{\bf z}

\newcommand{\rk}{{\rm rk}\,}
\newcommand{\icis}{{ICIS}}
\newcommand{\ind}{{\rm ind}}
\newcommand{\sgn}{{\rm sgn}}
\newcommand{\Sing}{{\rm Sing}\,}

\newcommand{\indGSV}{{\rm ind}_{\rm GSV}}
\newcommand{\indPH}{{\rm ind}_{\rm PH}}
\newcommand{\indPHN}{{\rm ind}_{\rm PHN}\,}
\newcommand{\indrad}{{\rm ind}_{\rm rad}}
\newcommand{\indhom}{{\rm ind}_{\rm hom}}

\newcommand{\Mmn}{M_{m,n}}

\newcommand{\funktion}[5]{\begin{array}{rccl} {#1}: & {#2} &
\longrightarrow & {#3}\\ & {#4} & \longmapsto & {#5} \end{array}}

\title[Indices of vector fields and 1-forms]{Indices of vector fields and 1-forms}
\author{Wolfgang Ebeling and Sabir M.~Gusein-Zade}
\thanks{Partially supported by DFG. The work of the second author was partially supported by the
grant 20-01-00579 of the Russian Foundation for Basic Research.}
\address{Wolfgang Ebeling, Institut f\"ur Algebraische Geometrie, Leibniz Universit\"at Hannover, Postfach 6009, D-30060 Hannover, Germany}
\email{ebeling@math.uni-hannover.de}
\address{Sabir M.~Gusein-Zade, Moscow State University, Faculty of mechanics and mathematics, Moscow Center for Fundamental and Applied Mathematics, GSP--1, Moscow, 119991, Russia \&
National Research University ``Higher School of Economics'',
Usacheva street 6, Moscow, 119048, Russia.}
\email{sabir@mccme.ru}

\date{}

\makeindex

\begin{document}

\maketitle

\begin{abstract} We discuss the notions of indices of vector fields and 1-forms and their generalizations to singular varieties and varieties with actions of finite groups, as well as indices of collections of vector fields and 1-forms.
\end{abstract}

\tableofcontents

\section*{Introduction}
Vector fields on a smooth manifold (a real or complex one) and their singular points play an important role in many different areas of mathematics. A classical invariant of a singular point of a vector field is its index.
The notion of the index  has a long history reflected in a number of classical sources. A famous result is the Poincar\'e--Hopf theorem which states that the sum of the indices of the (isolated) singular points of a vector field on a closed (compact, without boundary) manifold is equal to the Euler characteristic of the manifold. 

There is no straightforward 
generalization of the notion of the index of a singular point to vector fields on singular varieties. There are several concepts of indices in this case. Some of them require special conditions on the vector fields and/or on the singular variety.

In the smooth case, there is  essentially no difference between vector fields and 1-forms. In the case of singular varieties, these settings are essentially different. Traditionally the main attention was payed to indices of singular points of vector fields. A suggestion to consider indices of singular points of 1-forms alongside with (or instead of) indices of vector fields was first made by V.~I.~Arnold  in~\cite{Arnold79} (he used them for manifolds with boundaries; see also~\cite{Arnold04}). The authors started a study of indices of 1-forms on singular varieties \cite{EG01}.

Indices of vector fields or of 1-forms on (compact) complex analytic manifolds are related with the Euler characteristic, that is with the top Chern number. Other Chern numbers correspond to indices of collections of vector fields or of 1-forms. Therefore it is interesting to study such indices.

There are some equivariant versions of the Euler characteristic for spaces with an action of a finite group $G$. Therefore it is reasonable to try to define indices of singular points of $G$-invariant vector fields and 1-forms on $G$-varieties as well.

Here we give a survey on all these concepts. We give no proofs but for every statement we give precise references to articles where one can find more details including proofs. There is a comprehensive textbook \cite{BSS09} by J.-P.~Brasselet, J.~Seade, and T.~Suwa on a number of these subjects. We also wrote a survey article \cite{EGSchwerpunkt} about the developments till 2005. 

Let us outline the contents of this article. In the first section, we collect the basic notions and classical facts for the case of smooth manifolds.

In Section~\ref{sec:sing}, we discuss different generalizations of the notion of the index of a singular point to vector fields and 1-forms on singular varieties. 

We start (in Section~\ref{sect:rad}) with the notion which was classically defined first and which is very general. It is the radial or Schwartz index. The idea goes back to M.-H.~Schwartz who started
a comprehensive study of vector fields on singular (analytic) varieties  in \cite{Schwartz65}. She considered a class of vector fields important for a number of constructions related with vector fields on singular varieties: so called radial vector fields and radial extensions. For these vector fields she proved a version of the Poincar\'e--Hopf theorem \cite{Schwartz86a, Schwartz86b, Schwartz91}. Building on her work,
a general notion of an index of an isolated singular point of a vector field on a singular variety was introduced in~\cite{KT94}. However, this preprint was not published but only circulated around. (A revised version of it was only published almost 20 years later \cite{KT14}.)  Because of its restricted circulation, parts of it were later (re)elaborated in publications of other authors. In the general setting, for analytic varieties, this notion of index (called Schwartz index or radial index) was defined in \cite{EG99}, see also \cite{EGGD}. (For the case of varieties with isolated singularities, it was defined and studied, in particular, in \cite{ASV98,SS98}.)  It was noticed in \cite{Du10} that it can also be defined for semianalytic sets. The index was used to define characteristic classes for singular varieties, see, e.g., the surveys \cite{Brasselet00, Seade02}.

Section~\ref{sect:GSV} is devoted to a notion of an index which makes only sense on varieties of special types. It is the GSV index named after X.~G\'omez-Mont, J.~Seade, and A.~Verjovsky  who defined it in~\cite{GSV91} for vector fields on hypersurfaces with isolated singularities. It was generalized to vector fields on isolated complete intersection singularities in~\cite{SS96}. 

In Section~\ref{sect:PH}, we introduce the Poincar\'e--Hopf index for a vector field or 1-form on an isolated complex analytic singularity, a notion which is directly related to the Poincar\'e--Hopf theorem but is only defined for a smoothing of the singularity.

For the index of an isolated singular point of a holomorphic vector field on a complex manifold one has an algebraic formula, see Section~\ref{subsec:algebraic_form-smooth}. The search for such a formula for the index of an isolated singular point of a holomorphic vector field on a hypersurface with an isolated singularity was the motivation for G\'omez-Mont to introduce the notion of a homological index in \cite{GomezMont98}. This is the subject of Section~\ref{sect:Hom}.

The radial (Schwartz) index can be defined for singular points of vector fields on arbitrary analytic (or semi-analytic) varieties. Another notion of an index with this property is the Euler obstruction which is discussed in Section~\ref{sect:Eu}. The key ingredient in its construction was defined by R.~MacPherson in~\cite{MacPherson74}. He defined the Euler obstruction of the differential of the squared distance function. For vector fields, it was essentially defined by Brasselet and Schwartz in \cite{BS81}.

The homological index gives rise to algebraic formulas for the GSV index of a holomorphic vector field or 1-form on certain singular varieties generalizing those of Section~\ref{subsec:algebraic_form-smooth}. These formulas are discussed in Section~\ref{subsec:algebraic_form}. We also discuss analytic and topological formulas for the index.

Some of the indices are not defined for general analytic varieties, but  all the indices above are at least defined for isolated complete intersection singularities. The next more general class is the class of essentially isolated determinantal singularities introduced in \cite{EGSteklov} which recently attracted some attention. In the last subsection of Section~\ref{sec:sing} (Section~\ref{sect:det}), we study results on indices of 1-forms on such singularities.

In Section~\ref{sec:coll}, we consider analogues of the above indices for collections of vector fields and 1-forms. 
An analogue of the GSV index for them was introduced in \cite{EGBLMS}. It is discussed in Section~\ref{sect:GSVcoll}. 

An analogue of the notion of the Euler obstruction for collections of 1-forms corresponding to different Chern numbers leads to the notion of Chern obstructions introduced in \cite{EGBrasselet, EGTrieste}. This is considered in Section~\ref{sect:Chern}. There we also discuss relations between the Euler obstruction of a map and the Chern obstruction of a convenient collection of 1-forms observed by Brasselet, N.~G.~Grulha Jr., and M.~A.~S.~Ruas in \cite{BGR10}.

Finally we discuss, in Section~\ref{sect:Homcoll}, the generalization of the homological index to collections of 1-forms due to E.~Gorsky and the second author.

The final section (Section~\ref{sec:equiv_ind}) is devoted to the case that the variety carries an action of a finite group $G$.
Through the Poincar\'e--Hopf theorem indices of singular points of vector fields or 1-forms are often related with the Euler characteristic of the underlying variety. The Euler characteristic (properly defined) is an additive topological invariant of spaces of some kind, say, of locally closed unions of cells in finite CW-complexes. For topological spaces with additional structures one has other additive topological invariants which can be considered as generalized Euler characteristics. One can expect appropriate notions of indices of singular points corresponding to these concepts.
There are some notions of generalized Euler characteristics for spaces with an action of a  finite group $G$. Some of them take values in the group $\ZZ$ of integers: e.g.\ the alternating sum of the ranks of the invariant parts of the cohomology groups with compact support, the orbifold Euler characteristic (\cite{AS, HH}), etc. Another one is the alternating sum of classes of the cohomology groups as $G$-modules. It take values in the ring of representations of the group $G$. It was introduced in~\cite{Verdier}, see also~\cite{CTCWall}. The most general (in some sense~--- universal) concept of the generalized Euler characteristic for $G$-spaces is the equivariant Euler characteristic which takes values in the Burnside ring of the group $G$. It was introduced in~\cite{TtD}. We  discuss the definition in Section~\ref{subsec:equiv_Euler}.

A study of indices of singular points of $G$-invariant vector fields related with the Burnside ring of $G$ was started in~\cite{Luck}. The indices therein had values not in the Burnside ring of $G$, but in a related 
abelian group (not a ring) depending on the underlying $G$-manifold. A study of indices as elements of the Burnside ring itself was initiated in~\cite{EGEJM} and continued in~\cite{EGBrasil}. We first introduce these indices in the case of manifolds in Section~\ref{subsec:equiv_on_smooth}. 

The following subsections treat equivariant versions of the indices of vector fields and 1-forms on singular varieties discussed in Section~2. Namely, the equivariant radial index is treated in Section~\ref{subsec:equiv_radial_on_singular}, the equivariant GSV and Poincar\'e--Hopf index in Section~\ref{subsec:equiv_GSV}, the equivariant homological index in Section~\ref{subsec:equiv_homological} and the equivariant Euler obstruction in Section~\ref{subsec:equiv_Euler_obstruction}.

The final subsection Section~\ref{subsec:real_quotient}  discusses an attempt to generalize the Eisenbud--Levine--Khimshiashvili theorem to real quotient singularities.

\section{The case of smooth manifolds}\label{sec:smooth}
\subsection{The real index}\label{subsec:smooth_real}
Let $M$ be a smooth manifold of dimension $n$ and let $X$ be a vector field on $M$. A neighbourhood of a point $p \in M$ can be identified with a neighbourhood $U$ of the origin in $\RR^n$. In local coordinates around a point $p$, $X$ can be written as $X= \sum_{i=1}^n X_i(x) \frac{\partial}{\partial x_i}$, where $p$ corresponds to the origin in $\RR^n$. The vector field is called {\em continuous, smooth, analytic, etc.,}\index{vector field!continuous}\index{vector field!smooth}\index{vector field!analytic} if the functions $X_i$ are continuous, smooth, analytic, etc., respectively. A point $p \in M$ with $X(p)=0$ is called a {\em zero}\index{zero! of a vector field} or a {\em singular point}\index{singular point! of a vector field} of $X$. We shall define the index of a vector field at an isolated singular point.

For this purpose, we define the local degree of a mapping. Let $U \subset \RR^n$ be an open subset and 
$F: U \to \RR^n$ be a continuous mapping. Let $p \in U$ with $F(p)=0$ and let $B_\eps^n(p)=\{ x \in \RR^n \, | \, ||x-p||\leq \eps \}$ be the ball of radius $\eps$ centred at $p$  contained in $U$ such that there are no preimages of $F$ of the origin except $p$ inside $B_\eps^n(p)$. The {\em local degree} $\deg_p F$ \index{local degree of a mapping} of the mapping $F$ at the point $p$ is the degree of the mapping 
\[
\frac{F}{||F||}: S_\eps^{n-1}(p) \to S_1^{n-1}
\]
where $S_\eps^{n-1}(p)= \partial B_\eps^n(p)$ and $S_1^{n-1}$ is the unit sphere in $\RR^n$.

\begin{definition}
Let $X$ be a continuous vector field defined on $U \subset \RR^n$ and let $p \in U$ be an isolated zero of $X$. 
The {\em index} $\ind(X;\RR^n,p)$\index{index! of a real vector field} of $X$ at the singular point $p \in U$ is the degree of the mapping $(X_1, \ldots, X_n): U \to \RR^n$.
\end{definition}
One can easily see that the definition of the index is independent of the choice of the local coordinates and of the spheres.

Now let the vector field $X$ be smooth and $p$ a non-degenerate singular point of $X$. This means that $J_{X,p}:= \det \left( \frac{\partial X_i}{\partial x_j}(0) \right) \neq 0$. Then $\ind(X;\RR^n,p) = {\rm sign}\, J_{X,p}$, where ${\rm sign}\, J_{X,p}$ denotes the sign of $J_{X,p}$, i.e., 
\[ {\rm sign}\, J_{X,p}=\left\{ \begin{array}{cl} 1 &  \mbox{if }J_{X,p}>0, \\ -1 & \mbox{if } J_{X,p}<0. \end{array} \right.
\]
The index of an arbitrary isolated singular point $p$ of a smooth vector field $X$ is equal to the number of non-degenerate singular points $\widetilde{p}$ which split from the point $p$ under a generic perturbation $\widetilde{X}$ of the vector field $X$ in a neighbourhood of the point $p$ counted with the appropriate signs ${\rm sign}\, J_{\widetilde{X},\widetilde{p}}$.

One of the most important properties of the index of a vector field
is the Poincar\'e--Hopf theorem. Suppose that the manifold $M$ is
closed, i.e. compact without boundary, and that the vector field $X$
has finitely many singular points on it. This is equivalent to say that $X$ has only isolated singular points.

\begin{theorem}[Poincar\'e--Hopf]\label{theo1}\index{Poincar\'e--Hopf theorem}\index{theorem!Poincar\'e--Hopf}
Let $M$ be a closed (i.e.\ compact without boundary) manifold and $X$ be a smooth vector field with only finitely many singular points on it.
Then the sum
$$
\sum\limits_{p\in {\rm Sing}\, X}\ind(X;M,p)
$$
of indices of singular points of the vector field $X$ is equal to
the Euler characteristic $\chi(M)$ of the manifold $M$.
\end{theorem}
For a proof of this theorem see, e.g. \cite{Milnor65}. 

Instead of a vector field on $M$, one can consider a 1-form $\omega$ on $M$. In local coordinates, $\omega$ can be written as $\omega=\sum_{i=1}^n A_i(x) dx_i$. The notions {\em continuous, smooth, analytic, etc.,}\index{1-form!continuous}\index{1-form!smooth}\index{1-form!analytic} and {\em zero, singular point}\index{zero! of a 1-form}\index{singular point! of a 1-form} are defined analogously. Indeed, using a Riemannian metric one can identify vector fields and 1-forms on a smooth ($C^\infty$) manifold. In particular, the index ${\ind}_p\, \omega$\index{index! of a real 1-form} of $\omega$ at an isolated singular point $p$ is defined to be the degree of the 
map
$(A_1, \ldots , A_n):U \to \RR^n$.

\subsection{The complex index}\label{subsec:smooth_complex}
Now let $M$ be a complex manifold of complex dimension $n$ and $X$ be a vector field on $M$. The index $\ind(X;M,p)$\index{index! of a complex vector field} of the vector field $X$ at a singular point $p \in M$ is defined as the index of $X$ at $p$ on the underlying real $2n$-dimensional manifold. If the vector field is holomorphic and the singular point $p$ is non-degenerate, then the index is equal to $+1$. The index of an isolated singular point $p$ of a holomorphic vector field $X$ is positive. It is equal to the number of non-degenerate singular points which split from the point $p$ under a generic perturbation of the vector field in the neighbourhood of $p$.

If $M$ is a compact complex manifold, then its Euler characteristic $\chi(M)$ is equal to the characteristic number $\langle c_n(TM), [M] \rangle$, where $c_n(M)$ is the top Chern class of the manifold $M$. Therefore the Poincar\'e-Hopf theorem for a
 vector field $X$ on a compact complex manifold $M$ states that the sum of indices of the singular points on the vector field $X$ is equal to $\langle c_n(TM),[M] \rangle$.

Now let $\omega$ be a complex continuous 1-form on $M$. There is a one-to-one correspondence between complex 1-forms and real 1-forms on the underlying real $2n$-dimensional manifold. Namely, to a complex 1-form $\omega$ one can associate the real 1-form $\eta={\rm Re}\, \omega$. However, there is a difference in the orientation of the complex cotangent bundle $T^\ast M$ and the orientation of the real cotangent bundle with its complex structure forgotten. Therefore, the index $\ind({\rm Re}\, \omega; \RR^n,p)$ does not coincide with the index of the 1-form $\omega$ at the point $p$, but differs from it by the sign $(-1)^n$ (see, e.g., \cite[p.~235]{EGGD}). Therefore we define:

\begin{definition}
The {\em index} $\ind(\omega; \CC^n,p)$\index{index! of a complex 1-form} of the complex 1-form $\omega$ at a singular point $p$ is $(-1)^n$ times the index of the real 1-form ${\rm Re}\, \omega$ at $p$:
\[
\ind(\omega; \CC^n,p) := (-1)^n \ind({\rm Re}\, \omega; \RR^n,p).
\]
\end{definition}

With this definition, the Poincar\'e-Hopf theorem for a complex 1-form $\omega$ on a compact complex manifold $M$ states that the sum of the indices of the singular points of $\omega$ is equal to $(-1)^n \chi(M)=\langle c_n(T^\ast M), [M] \rangle$.

\subsection{Collections of sections of a vector bundle} \label{sect:coll}
The complex index introduced above is connected with
the Euler characteristic, hence with the characteristic number $\langle c_n(TM), [M] \rangle$, where $c_n(M)$ is the top Chern class of the manifold $M$. In this section, we discuss indices related to other characteristic numbers.

Let $M$ be a complex analytic manifold of dimension $n$ and let $\{ X_j \}=(X_1, \ldots , X_k)$ be a collection of $k$ continuous vector fields on $M$. A {\em singular point}\index{singular point! of a collection} of $\{ X_j \}$ is a point $p \in M$ where $(X_1(p), \ldots , X_k(p))$ are linearly dependent. A  collection $\{ X_j \}$ is also called a {\em $k$-field}\index{$k$-field} and a non-singular one is called a {\em $k$-frame}\index{$k$-frame}. We recall the construction of Chern classes by obstruction theory.

For natural numbers $n$ and $k$ with $n \geq k$, let $M_{n,k}$ be
the space of $n \times k$ matrices with complex entries and let $D_{n,k}$
be the subspace of $M_{n,k}$ consisting of matrices of rank less than
$k$. The subset $D_{n,k}$ is a subvariety of $M_{n,k}$
of codimension $n-k+1$. The complement $W_{n,k} = M_{n,k}
\setminus D_{n,k}$ is the Stiefel manifold of $k$-frames (collections of $k$ linearly independent vectors) in $\CC^n$. It is known that $W_{n,k}$
is $(2n-2k)$-connected and $H_{2n-2k+1}(W_{n,k}) \cong \ZZ$ (see, e.g.,
\cite{Husemoller75}). The latter fact also implies that the subvariety $D_{n,k}$ is irreducible. Since
$W_{n,k}$ is the complement of an irreducible complex analytic subvariety
of codimension $n-k+1$ in $M_{n,k}$, there is a natural choice of a
generator of the homology group $H_{2n-2k+1}(W_{n,k}) \cong \ZZ$. Namely,
the (``positive'') generator is the boundary of a small ball  in a smooth complex analytic slice in $M_{n,k} \cong \CC^{nk}$ transversal to the irreducible subvariety $D_{n,k}$ at a non-singular point (oriented in the standard way).

Let $\{ X_j \}$ be a $k$-field on $M$ with isolated singular points.
Let (K) be a triangulation of $M$ such that all singular points of $\{ X_j \}$ are vertices of (K) and let  (D) be a cell decomposition of $M$ dual to (K). Let $\sigma$ be a $2(n-k+1)$-cell of (D) which is contained in an open subset $U \subset M$ where the tangent bundle $TM$ is trivial. Let
\[ (X_1(y), \ldots , X_k(y) )
\]
be the $n \times (n-k+1)$-matrix the columns of which consist of the components of the vectors $X_1(y), \ldots , X_k(y)$ with respect to this trivialization. Let $\psi_\sigma : \partial \sigma \cong S^{2n-1} \to W_{n,n-k+1}$ be the mapping which sends a point $y \in \partial \sigma$ to the matrix $(X_1(y), \ldots , X_k(y) )$.

\begin{definition} The {\em index} $\ind(\{ X_j \};\sigma)$ of the $k$-field $\{ X_j \}$ on $\sigma$ is the degree of the map $\psi$, i.e., the obstruction to extend the $k$-frame $\{ X_j \}$ from the boundary $\partial \sigma$ of the cell $\sigma$ to its interior.
\end{definition}

This defines a cochain $\gamma \in C^{2(n-k+1)}(M;\ZZ)$ by setting $\gamma(\sigma)=\ind(\{ X_j \};\sigma)$ for each $2(n-k+1)$-cell and extending it linearly.
This cochain is in fact a cocycle and represents the Chern class $c^{n-k+1}(M)$ of $M$.

Let $\pi: E \to M$ be a complex
analytic vector bundle of rank $m$ over a complex analytic
manifold $M$ of dimension $n$. (Special cases of interest are the tangent and the cotangent bundles of $M$.) We shall now generalize the construction above.

Let $\{\omega^{(i)}_j\}$ ($i=1, \ldots , s$; $j=1, \ldots , m-k_i+1$;
$\sum_{i=1}^s k_i = n$) be a collection of continuous sections of the vector bundle
$\pi:E \to M$. A point $p \in M$ is called non-singular for the collection  $\{\omega^{(i)}_j\}$ if at least for some $i \in \{1, \ldots, s\}$
the values $\omega^{(i)}_1(p), \ldots , \omega^{(i)}_{m-k_i+1}(p)$ are linearly independent. This means that for this $i$ the vectors $\omega^{(i)}_1(p), \ldots , \omega^{(i)}_{m-k_i+1}(p)$ form an $(m-k_i+1)$-frame. We shall assume that the collection $\{\omega^{(i)}_j\}$ has only isolated singular points. We shall define an index for such a collection, cf.\ \cite{EGBLMS}.

Let ${\bf
k}=(k_1, \ldots , k_s)$ be a sequence of positive integers with
$\sum_{i=1}^s k_i = k$. Consider the space $M_{m, {\bf k}}=
\prod_{i=1}^s M_{m,m-k_i+1}$ and the subvariety $D_{m, {\bf k}}=
\prod_{i=1}^s D_{m,m-k_i+1}$ in it. The variety $D_{m, {\bf k}}$ consists
of sets $\{A_i\}$ of $m \times (m-k_i+1)$ matrices such that
$\mbox{rk}\, A_i < m-k_i+1$  for each $i=1, \ldots , s$. Since
$D_{m, {\bf k}}$ is irreducible of codimension $k$, its complement
$W_{m, {\bf k}}= M_{m, {\bf k}} \setminus D_{m, {\bf k}}$
is $(2k-2)$-connected, $H_{2k-1}(W_{m, {\bf k}}) \cong \ZZ$, and there
is a natural choice of a generator of the latter group. This choice
defines a degree (an integer) of a map from an oriented manifold of
dimension $2k-1$ to the manifold $W_{m, {\bf k}}$.

Let us choose a trivialization of the vector
bundle $\pi: E \to M$ in a neighbourhood of a point $p$, let
\[(\omega^{(i)}_1(x), \ldots, \omega^{(i)}_{m-k_i+1}(x))
\]
 be the
$m \times (m-k_i+1)$-matrix the columns of which consist of the components
of the sections $\omega^{(i)}_j(x)$, $j=1, \ldots , m-k_i+1$, $x \in M$,
with respect to this trivialization. Let $\Psi_p$ be the mapping from
a neighbourhood of the point $p$ to $M_{m, {\bf k}}$ which sends
a point $x$ to the collection of matrices $\{ (\omega^{(i)}_1(x), \ldots,
\omega^{(i)}_{m-k_i+1}(x))\}$, $i=1, \ldots, s$.
Its restriction $\psi_p$ to a small sphere $S^{2n-1}_\eps(p)$ around
the point $p$  maps this sphere to the subset $W_{m, {\bf k}}$. The sphere
$S^{2n-1}_{\eps}(p)$ is oriented as the boundary of the corresponding
ball in the complex affine space $\CC^n$.

\begin{definition} The {\em
index} $\ind(\{\omega^{(i)}_j\};M,p)$\index{index! of a collection}
of the collection of sections $\{\omega^{(i)}_j\}$ at the point $p$ is the degree of the mapping
$\psi_p: S^{2n-1}_\eps(p) \to W_{m, {\bf k}}$.
\end{definition}

One can deform the collection $\{\omega^{(i)}_j\}$ of sections or,
equivalently, the map $\Psi_p$ so that this map becomes smooth and transversal
to the variety $D_{m, {\bf k}}$ (at smooth points of the latter one).
This implies that
the index $\ind(\{\omega^{(i)}_j\};M,p)$ is equal
to the intersection number of the germ of the image of the map $\Psi_p$ with the variety $D_{m, {\bf k}}$.

\begin{definition} A singular point $p$ of the collection $\{\omega^{(i)}_j\}$ of sections is called
{\em non-degenerate}  \index{singular point!non-degenerate} if the map $\Psi_p$
is  smooth and transversal to the variety
$D_{n, {\bf k}} \subset {\mathcal M}_{n, {\bf k}}$ at a non-singular point of it.
\end{definition}

If $p$ is a non-degenerate singular point of the collection $\{\omega^{(i)}_j\}$, then $\mbox{ind}_p\{\omega^{(i)}_j\}=\pm 1$. If all the sections $\omega^{(i)}_j$ are complex analytic, then this index is equal to $+1$.

The following statement is a generalization of the well known fact that the
($2(n-k)$-dimensional) cycle Poincar\'e dual to the characteristic
class $c_k(E)$ ($k=1, \ldots , m$) is represented by the set of
points of the manifold $M$ where $m-k+1$ generic sections of the
vector bundle $E$ are linearly dependent (cf., e.g., \cite[p.~413]{GH78}).

\begin{theorem} \label{thm:Chernsec}
Let $\sum_{i=1}^s k_i  = n$ and suppose that the collection
$\{\omega^{(i)}_j\}$ $($$i=1, \ldots , s$; $j=1, \ldots, m-k_i+1$$)$ of
sections of the vector bundle $\pi: E \to M$ over a closed complex
manifold $M$  has only isolated singular points. Then the sum of the
indices of these points is equal to the characteristic number
$\langle\prod_{i=1}^s c_{k_i}(E), [M]\rangle$ of the vector bundle $E$.
\end{theorem}

\subsection{Algebraic formulas for the indices}\label{subsec:algebraic_form-smooth}
It turns out that, if the vector fields (or 1-forms) under consideration are analytic (and, in the real case, the singular points are algebraically isolated) one has algebraic formulas for the indices considered above. 

The simplest algebraic formula is for the index of an isolated singular point of a holomorphic vector field on a complex manifold. In local coordinates $z=(z_1, \ldots , z_n)$ centred at the singular point, a vector field can be written as $X=\sum_{i=1}^n X_i(z) \frac{\partial}{\partial z_i}$, where the function germs $X_i$ are holomorphic. Let $\calO_{\CC^n,0}$ be the ring of germs of holomorphic functions of $n$ variables. 

The first proof of the following theorem is usually attributed to Palamodov \cite{Palamodov67}. Without a detailed proof, it was known earlier.

\begin{theorem}[Palamodov]
The index $\ind(X; \CC^n,0)$ of the singular point of the holomorphic vector field $X$ is equal to the dimension of the complex vector space $\calO_{\CC^n,0}/(X_1, \ldots , X_n)$, where $(X_1, \ldots , X_n)$ is the ideal generated by the germs $X_1, \ldots , X_n$.
\end{theorem}

Recall that $J_{X,0}=\det \left( \frac{\partial X_i}{\partial z_j}(0) \right)$ denotes the determinant of the Jacobian matrix of  $(X_1, \ldots , X_n)$.
Then one has the following residue formula for the index
\[ \ind(X; \CC^n,0) = {\rm Res}\left[ 
\begin{array}{c} J_{X,0} d{\zz} \\
X_1 \cdots X_n \end{array} \right] =  \frac{1}{(2 \pi i)^n} \int_\Gamma \frac{J_{X,0}}{X_1 \cdots X_n} d{\zz}
\]
where $d{\zz} :=dz_1 \wedge \cdots \wedge d z_n$, $\Gamma$ is the real $n$-cycle $\{\Vert X_k \Vert=\delta_k, k=1, \ldots ,n\}$ for positive $\delta_k$ small enough, and $\Gamma$ is oriented so that $d(\arg X_1) \wedge \cdots \wedge d(\arg X_n) \geq 0$, see also \cite{BB70}.

For a real analytic vector field such that its complexification has an isolated singular point, the index can be computed as the signature of a certain quadratic form: \cite{EL77, Kh77}.

Let $F=(f_1, \ldots , f_n): (\RR^n,0) \to (\RR^n,0)$ be the germ of an analytic mapping such that
$F_\CC^{-1}(0)=0$, 
where $F_\CC: (\CC^n,0) \to (\CC^n,0)$ is the complexification of $F$. (In this situation one says that $F$ has an algebraically isolated\index{isolated!algebraically} preimage of the origin.) Let $\calE_{\RR^n,0}$ be the ring of germs of analytic functions on $(\RR^n,0)$. By the assumption  
$F_\CC^{-1}(0)=0$,  
the factor algebra $Q_F:=\calE_{\RR^n,0}/(f_1, \ldots , f_n)$ has finite dimension. (This dimension is equal to $\dim \calO_{\CC^n,0}/(f_1, \ldots , f_n)=\deg_0F_\CC$.) We consider on $Q_F$ the natural residue pairing
\[
\funktion{B_F}{Q_F \times Q_F}{\RR}{(\varphi,\psi)}{{\rm Res}\left[ 
\begin{array}{c} \varphi({\xx}) \psi({\xx}) d{\xx} \\
f_1 \cdots f_n \end{array} \right]} 
\] 
and the residue is similarly defined as in the complex case, namely 
\[
{\rm Res}\left[ 
\begin{array}{c} \varphi({\xx}) \psi({\xx}) d{\xx} \\
f_1 \cdots f_n \end{array} \right] 
 =  \frac{1}{(2 \pi i)^n} \int \frac{\varphi({\xx}) \psi({\xx})}{f_1 \cdots f_n} d{\xx}
\]
where $d{\xx}:=dx_1 \wedge \cdots \wedge dx_n$ and the integration is along
the cycle given by the equations $\Vert f_k(\xx)\Vert=\delta_k$ with positive $\delta_k$ small enough.

\begin{theorem}[Eisenbud--Levine--Khimshiashvili]\label{theo3}\index{Eisenbud--Levine--Khimshiashvili theorem}\index{theorem!Eisenbud--Levine--Khimshiashvili}
The degree $\deg_0 F$ of the map germ $F$
is equal to the signature $\sgn \, B_F$ of the quadratic form $B_F$.
\end{theorem}

For a proof of this theorem see also \cite{AGV85}.

This can also be
interpreted as a formula for the index of the singular point of the vector field
$X:=\sum f_i \frac{\partial}{\partial x_i}$ or of the 1-form $\omega:=\sum f_i dx_i$. Moreover, the choice
of a volume form permits to identify the algebra $Q_F$ (as a vector space) with the space
$\Omega_\omega= \Omega^n_{\RR^n,0}/\omega \wedge \Omega^{n-1}_{\RR^n,0}$. 

Now let $\{\omega^{(i)}_j\}$ ($i=1, \ldots , s$; $j=1, \ldots , m-k_i+1$;
$\sum_{i=1}^s k_i = n$) be a collection of holomorphic sections of the (trivial) vector bundle
$\pi:\CC^m\times(\CC^n,0) \to (\CC^n,0)$ with an isolated singular point at the origin (see Subsection~\ref{sect:coll}).
Let $I_{\{\omega^{(i)}_j\}}$ be the
ideal in the ring ${\mathcal O}_{\CC^n,0}$ of germs of analytic functions
of $n$ variables generated by the $(m-k_i+1) \times (m-k_i+1)$-minors
of the matrices $(\omega^{(i)}_1, \ldots, \omega^{(i)}_{m-k_i+1})$ for
all $i=1, \ldots, s$. Then one has the following algebraic formula for the index (see \cite[Theorem~2]{EGBLMS}).

\begin{theorem} \label{Th2}
One has
$$
\ind(\{\omega^{(i)}_j\};\CC^n,0) =
\dim_\CC {\mathcal O}_{\CC^n,0}/I_{\{\omega^{(i)}_j\}}.
$$
\end{theorem}

\section{Vector fields and 1-forms on singular varieties}\label{sec:sing}
\subsection{Radial index} \label{sect:rad}
Let $V$ be  a closed (real) subanalytic variety in  a smooth manifold $M$, where $M$ is equipped with a (smooth) Riemannian metric.
Let $V=\bigcup_{i=1}^q V_i$ be a subanalytic Whitney stratification of $V$ (see \cite{Trotman20} for this notion).
A (continuous) {\em stratified vector field}\index{vector field!stratified} on $V=\bigcup_{i=1}^q V_i$
is a vector field such that, at each point $p$ of $V$, it is tangent
to the stratum containing $p$.

\begin{definition}
The germ $X$ of a vector field on the germ $(V,p)$ is called {\em radial}\index{radial vector field}\index{vector field! radial} if, for all $\eps > 0$ small enough, the vector field is transversal to the boundary of the $\eps$-neighbourhood of the point $p$ and is directed outwards.
\end{definition}


Let $p\in V$, let $V_{(p)}=V_i$ be the stratum containing $p$, $\dim V_i=k$,
and let $X$ be a stratified vector field on $V$ in a neighbourhood
of the point $p$. The following definition was given by Schwartz in  \cite{Schwartz86a}.

\begin{definition}
The vector field $X$ is called a {\em radial extension}\index{radial extension! of a vector field} of the vector field
$X_{\vert V_i}$ if, for all $\eps >0$ small enough, it is transversal to the boundary of the $\eps$-tubular neighbourhood of $V_i$ and points outwards of the neighbourhood.
\end{definition}

\begin{remark}
Note that the definitions of the radial extension in \cite[Definition~2.3.2]{BSS09}), \cite[p.~144]{EG99}, and \cite[p.~288]{EGEJM} use different formulations, but are somewhat inaccurate.
\end{remark}

The existence of a radial extension of a vector field on a stratum $V_i$  is proved in \cite[Section~III.7]{BS81} (see also \cite[Lemme~3.1.2]{Schwartz91}).

Let $X$ be a stratified vector field on $(V,p)$ with an isolated singular point
(zero) at the origin. Let $B_\eps$ be a closed $\eps$-neighbourhood in the ambient Riemannian manifold $M$ around the point $p$ (small enough so that the boundary $\partial B_\eta$ of the $\eta$-neighbourhood of $p$ with $0< \eta \leq \eps$ intersects all the strata $V_i$ transversally and the vector field $X$ has no singular points on $(V \setminus \{ p \}) \cap B_\eps$ ). One can show that there exists a (continuous)
stratified vector field $\widetilde{X}$ on $V$ such that:
\begin{itemize}
\item[(1)] the vector field $\widetilde{X}$ coincides with $X$ on a neighbourhood
of the intersection of $V$ with the boundary $\partial B_\eps$ of the $\eps$-neighbourhood 
 around the point $p$;
\item[(2)] the vector field has a finite number of singular points (zeros);
\item[(3)] in a neighbourhood of each singular point $p\in V\cap B_\eps$, $p\in V_i$,
the vector field $\widetilde{X}$ is a radial extension of
its restriction to the stratum $V_i$.
\end{itemize}

\begin{definition}
The {\em radial index} (or {\em Schwartz index})\index{radial index! of a vector field}\index{Schwartz index}
$\indrad(X;V,p)$ of the vector field $X$
on $V$ at the point $p$ is 
$$
\indrad(X;V,p):=
\sum_{{p}\in {\rm Sing}\widetilde{X}} 
{\rm ind\,}(\widetilde{X}_{\vert V_{(p)}};V_{(p)},p)\,,
$$
where ${\rm ind\,}(\widetilde{X}_{\vert V_{(p)}},V_{(p)},p)$ is the usual index of the 
restriction of the vector field $\widetilde{X}$ to the smooth manifold $V_{(p)}$.
\end{definition}

The described notions of  a radial vector field and a radial extension depend on the choice of the Riemannian metric, but the radial index does not depend on this choice, because, for instance, the space of Riemannian metrics is pathwise connected.

Now let $\omega$ be (the germ  at $p$ of) a (continuous) 1-form on $(V,p)$, i.e.\ the restriction to $V$ of a 1-form defined in a neighbourhood of the point $p$ in the ambient manifold $M$. Let $V=\bigcup_{i=1}^q V_i$ be a subanalytic Whitney stratification of $V$. A point $p \in V$ is a {\em singular point} of $\omega$ if the restriction of $\omega$ to the stratum $V_{(p)}$ containing $p$ vanishes at the point $p$.

\begin{definition}
The germ $\omega$ of a 1-form at the point $p$ is called {\em radial}\index{radial 1-form}\index{1-form!radial} if, for all $\eps$ small enough, the 1-form is positive on the outward normals to the boundary of the $\eps$-neighbourhood of the point $p$.
\end{definition}

An example of a radial 1-form is the germ of the 1-form $d\rho^2$, where $\rho$ is the distance function from $p$ induced by the Riemannian metric.

\begin{remark} Note that the initial definition of a radial 1-form in \cite{EGGD} used other words, but was somewhat inaccurate. This was noticed by an anonymous referee of the paper \cite{Gusein-Zade21}.
\end{remark}

Let $p$, $p \in V_i=V_{(p)}$, $\dim V_{(p)}=k$, and let $\eta$ be a 1-form defined in a neighbourhood of the point $p$. As above, let $N_i$ be a normal slice (with respect to the Riemannian metric) of $M$ to the
stratum $V_i$ at the point $p$ and $h$ a diffeomorphism
from a neighbourhood of $p$ in $M$ to the product $U_i(p)\times N_i$, where
$U_i(p)$ is an $\eps$-neighbourhood of $p$ in $V_i$, which is the identity on $U_i(p)$.

\begin{definition}
A 1-form $\eta$ is called a {\em radial extension}\index{radial extension! of a 1-form} of the 1-form $\eta_{\vert_{V_{(p)}}}$ if there exists such a diffeomorphism $h$ which identifies $\eta$ with the restriction to $V$ of the 1-form $\pi_1^\ast \eta_{\vert_{V_{(p)}}} + \pi_2^\ast \eta^{\rm rad}_{N_i}$, where $\pi_1$ and $\pi_2$ are the projections from a neighbourhood of $p$ in $M$ to $V_{(p)}$ and $N_i$ respectively and $\eta^{\rm rad}_{N_i}$ is a radial 1-form on $N_i$.
\end{definition}

For a 1-form $\omega$ on $(V,p)$ with an isolated singular point at the point $p$ there exists a 1-form $\widetilde{\omega}$ on $V$ which possesses the obvious analogues of the properties (1)--(3) of the vector field $\widetilde{X}$ above.

\begin{definition}
The {\em radial index}\index{radial index! of a 1-form} $\indrad(\omega;V,p)$ of the 1-form $\omega$ at the point $p$ is 
$$
\indrad(\omega;V,p)=
\sum_{{p}\in {\Sing}\widetilde{\omega}} 
\ind(\widetilde{\omega}_{\vert V_{(p)}};V_{(p)},p)\,,
$$
where $\ind(\widetilde{\omega}_{\vert V_{(p)}};V_{(p)},p)$ is the usual index of the restriction of the 1-form $\widetilde{\omega}$ to the stratum $V_{(p)}$.
\end{definition}

The definition of the radial index does not depend on the stratification and on the chosen vector field $\widetilde{X}$ nor of the chosen 1-form $\widetilde{\omega}$. The sum of the indices of an appropriate deformation of the vector field or 1-form on the strata is the same for a stratification and for a refinement of it. This implies that the radial index does not depend on the stratification. (For a vector field one can consider the intersection of two stratifications. For a 1-form one can consider the minimal Whitney stratification of the variety.) Moreover, the radial index does not depend on the chosen deformation of the vector field or 1-form. This follows from the following proposition which is proved in \cite[Proposition~2.1]{EGEJM}.

\begin{proposition}\label{prop-distr}
The number of singular points (counted with multiplicities) of the vector field $\widetilde{X}$ or of a 1-form $\widetilde{\omega}$
on a fixed stratum $V_i$ does not depend on the choice of the vector field $\widetilde{X}$ or of the 1-form $\widetilde{\omega}$ respectively (and therefore only depends on $X$ or $\omega$ respectively).
\end{proposition}

Therefore the radial index is well-defined.

It follows from the definition that the radial index satisfies the {\em law of conservation of number}. For a vector field $X$ this means the following: if a vector field $X'$ with isolated singular points on $V$ is close to the vector field $X$, then
\[
\indrad(X;V,p) = \sum_{q\in \Sing X'} \indrad(X';V,q),
\]
where the sum on the right hand side runs over all singular points $q$ of the vector field $X'$ on $V$ in a neighbourhood of $p$.

The radial index generalizes the usual index for vector fields or 1-forms on a smooth manifold. In particular, one has a generalization of the Poincar\'e-Hopf theorem:

\begin{theorem}[Poincar\'e--Hopf]
For a compact real subanalytic variety $V$ and a vector field $X$ or a 1-form $\omega$ with isolated singular points on $V$, one has
$$
\sum_{Q \in \Sing X} {\ind}_{\rm rad}\,(X;V,Q) =
\sum_{Q \in \Sing \omega} {\ind}_{\rm rad}\,(\omega;V,Q) = \chi(V)
$$
where $\chi(V)$ denotes the Euler characteristic of the space (variety) $V$.
\end{theorem}

In \cite{Lap18}, the radial index of a vector field with an isolated zero on a real closed semialgebraic set with an isolated singularity is related to an intersection index.

Now we consider the germ $(V,p)$ of a complex analytic variety of pure dimension $n$ embedded in the germ of a complex manifold $(M,p)$. Then one can define analogously the notion of a radial index for a singular point of a complex vector field or complex 1-form on $(V,p)$. In particular, for the germ of a complex 1-form $\omega$ on $(V,0)$ the radial index $\indrad(\omega;V,P)$ is $(-1)^n$ times the radial index $\indrad({\rm Re}\,\omega;V,P)$ of the real 1-form ${\rm Re}\,\omega$ on $(V,P)$.

GSV 
\subsection{GSV index} \label{sect:GSV}
Let $(V,0)\subset(\CC^{N},0)$ 
be the germ of an $n$-dimensional  complete intersection with an isolated singularity at the origin, defined by a holomorphic map germ
\[
f=(f_1, \ldots , f_{N-n}): (\CC^{N},0) \to (\CC^{N-n},0),
\]
i.e.\ $V=f^{-1}(0)$. The germ $(V,0)$ is called
an isolated complete intersection singularity (abbreviated {\icis}\ in the sequel). Let $z_1, \ldots , z_{N}$ denote the coordinates of $\CC^{N}$. Let $X=\sum_{i=1}^{N} X_i(z) \frac{\partial}{\partial z_i}$ be the germ of a (continuous) vector field on $(\CC^{N},0)$ tangent to $V$, i.e.\ $X(z) \in T_zV$ for all $z \in V \setminus \{ 0 \}$. Suppose that $X$ has an isolated singular point at the origin. Then the following index is defined. It was introduced for hypersurface singularities by X.~G\'omez-Mont, J.~Seade, and A.~Verjovsky \cite{GSV91} and generalized to  {\icis}\ by Seade and T.~Suwa \cite{SS98}. It is called GSV index.

Let $B_\eps \subset \CC^{N}$ denote the ball of radius $\eps$ centred at the origin. Let $\eps > 0$ be chosen small enough so that the functions $f_1, \ldots f_{N-n}$ and the vector field $X$ are defined in a neighbourhood of $B_\eps$, $V$ is transversal to the sphere $S_\eta= \partial B_\eta$ for $0 < \eta \leq \eps$, and the vector field $X$ has no zeros on $V$ inside the ball $B_\eps$ except possibly the origin. The intersection $K:=V \cap S_\eps$ is called the {\em link} of the ICIS $(V,0)$. The link $K$ is a $(2n-1)$-dimensional manifold and has a natural orientation as the boundary of the complex manifold $(V \cap B_\eps) \setminus \{ 0 \}$. 

Define the {\em gradient vector field}\index{gradient vector field}\index{vector field!gradient} ${\rm grad}\, f_i$ of a function germ $f_i$ by
\[
{\rm grad}\, f_i= \left(\,\overline{\frac{\partial f_i}{\partial z_1}}, \ldots ,
\overline{\frac{\partial f_i}{\partial z_{N}}}\,\right).
\]
Note that it depends on the choice of the coordinates $z_1, \ldots , z_{N}$. The gradient vector fields ${\rm grad}\, f_1, \ldots , {\rm grad}\, f_{N-n}$ are linearly  independent everywhere on $V$ except (possibly) at the origin. The set $\{ X(z), {\rm grad}\, f_1(z), \ldots , {\rm grad}\, f_{N-n}(z) \}$ is an $(N-n+1)$-frame  at each point of $K$. This frame defines a continuous map
\begin{equation*}\label{eq_psi}
\Psi=(X,{\mbox{grad\,}} f_1, \ldots , {\mbox{grad\,}} f_{N-n}):
K\to W_{N, N-n+1}
\end{equation*}
from the link $K$ to the Stiefel manifold $W_{N, N-n+1}$ of complex  $(N-n+1)$-frames in $\CC^{N}$. It is known that the Stiefel manifold $W_{N, N-n+1}$ is $2(n-1)$-connected and $H_{2n-1}(W_{N, N-n+1}) \cong \ZZ$, see Section~\ref{sect:coll}. There is a natural choice of the generator of $H_{2n-1}(W_{N, N-n+1}) \cong \ZZ$.  Therefore we can make the following definition:

\begin{definition}
The {\em GSV index } \index{GSV index! of a vector field} $\indGSV(X;V,0)$ of the vector field
$X$
on the {\icis} $V$ at the origin is the degree of the map
$$
\Psi: K\to W_{N,N-n+1}\,.
$$
\end{definition}

\begin{remark}
Note that one uses the complex conjugation for this definition.
Therefore the components of the discussed map are of different tensor nature.
Whereas $X$ is a vector field, ${\mbox{grad\,}}f_i$ is more similar
to a covector.
\end{remark}

One can also consider the map $\Psi$ as a map from $V$ to the space $M_{N,N-n+1}$ of $N \times (N-n+1)$ matrices with complex entries (defined in a neighbourhood of the ball $B_\eps$). It maps the set $V \setminus \{ 0\}$ to the Stiefel manifold  $W_{N, N-n+1}=M_{N,N-n+1}\setminus D_{N,N-n+1}$ (see Section~\ref{sect:coll}). Therefore we have the following result:

\begin{proposition}\label{prop_gsv1}
The GSV index $\ind_{\rm GSV} \, (X;V,0)$ of the vector field $X$
on the {\icis} $V$ at the origin is equal to the intersection number
$(\Psi(V)\circ D_{N, N-n+1})$ of the image $\Psi(V)$ of the {\icis}
$V$ under the map $\Psi$ and the variety $D_{N, N-n+1}$ at the origin.
\end{proposition}

Note that, even if the vector field $X$ is holomorphic, the image $\Psi(V)$  is in general not a complex analytic variety because we use the complex conjugation in the definition of $\Psi$.

Now we consider the case of a 1-form on $(V,0)$. Let $\omega=\sum A_i(z)dz_i$  be
a germ of a continuous 1-form on $(\CC^{N}, 0)$ which as a 1-form
on the {\icis} $V$ has (at most) an isolated singular point at the origin
(thus it does not vanish on the tangent space $T_pV$ to the variety
$V$ at all points $p$ from a punctured neighbourhood of the origin
in $V$). The set $\{ \omega(z), df_1(z), \dots , df_{N-n}(z) \}$ is a $(N-n+1)$-frame in the space dual to $\CC^{N}$ for all $z \in K$. Therefore one has a map
$$
\Psi=(\omega, df_1, \ldots , df_{N-n}) : K \to W_{N, N-n+1}.
$$
Here $W_{N, N-n+1}$ is the Stiefel manifold of $(N-n+1)$-frames in
the space dual to $\CC^{N}$.

\begin{definition}
The {\em GSV index} \index{GSV index! of a 1-form} $\ind_{\rm GSV} \, (\omega;V,0)$ of the 1-form $\omega$
on the {\icis} $V$ at the origin is the degree of the map
$$
\Psi: K\to W_{N,N-n+1}\,.
$$
\end{definition}

Just as above $\Psi$ can be considered as a map from the {\icis} $V$
to the space $M_{N, N-n+1}$ of $N \times (N-n+1)$-matrices and the GSV index of $\omega$ is equal to the intersection number $(\Psi(V)\circ D_{N, N-n+1})$. In contrast to the case of a vector field,
if the 1-form $\omega$ is holomorphic, the map $\Psi$ and the set
$\Psi(V)$ are complex analytic.

A similar construction can be considered in the real setting, see \cite{ASV98} for more details.

\subsection{Poincar\'e--Hopf index} \label{sect:PH}
Let $(V,0) \subset (\CC^N,0)$ be the germ of a purely $n$-dimensional complex analytic variety with an isolated singularity at the origin. A {\em smoothing} of $(V,0)$ is a 1-parameter deformation $F: (\calV,0) \to (\CC,0)$ of $(V,0)$ (that is $F^{-1}(0)=(V,0)$) such that for $t \in \CC \setminus \{ 0 \}$ sufficiently close to 0 the fibre $\calV_t=F^{-1}(t)$ is smooth, cf., e.g.\ \cite[Definition~7.3.1]{Greuel20}. The germ $(V,0)$ is called a {\em smoothable singularity} if there exists a smoothing of $(V,0)$. 

Let $(V,0)$ be a smoothable singularity and let $F: \calV \to \CC$ be a suitable representative of a smoothing of $(V,0)$. For simplicity, we assume that $\calV$ is embedded in an open neighbourhood $U \subset \CC^N$ of the origin.  Denote by $B_\eps$ the ball of radius $\eps$ centred at the origin in $\CC^{N}$ and by $\Delta_\eta$ the disc in $\CC$ of radius $\eta$ centred at 0. Let $\Delta_\eta^\ast=\Delta_\eta \setminus \{ 0 \}$. By \cite[Theorem~1.1]{Le76} (see also \cite[Theorem~6.4.1]{LNS20}), for $\eps \gg \eta > 0$ sufficiently small, the mapping 
\[ 
F_{\vert_{F^{-1}(\Delta_\eta^\ast) \cap B_\eps}} : F^{-1}(\Delta_\eta^\ast) \cap B_\eps \to \Delta_\eta^\ast
\]
is the projection of a differentiable fibre bundle over $\Delta_\eta^\ast$. Let $V^F_t:= \calV_t \cap B_\eps$ be the (Milnor) fibre of this bundle over $t \in \Delta_\eta^\ast$.

Now let $X$ be the germ of a continuous vector field on $(V,0)$ with an isolated singularity at $0$. Then the vector field does not vanish on $V \cap S_\eps$, where $S_\eps=\partial B_\eps$ is the boundary of the ball $B_\eps$. Moreover, the intersection $V \cap S_\eps$ is isotopic to the intersection of $V^F_t$ with this sphere. Therefore we can assume that the vector field $X$ is defined on the boundary $\partial V^F_t$ of the Milnor fibre. By \cite[Theorem~1.1.2]{BSS09}, there exists an extension $\widetilde{X}$ of the vector field  $X$ to the interior of the Milnor fibre $V^F_t$ with a finite number of singular points.

\begin{definition} The {\em Poincar\'e--Hopf index}\index{Poincar\'e--Hopf index} of $X$ on $(V,0)$ relative to the smoothing $F$ is
\[
\indPH^F(X;V,0) := \sum_{q \in \Sing \widetilde{X}}  \ind(\widetilde{X};V^F_t,q),
\]
where the sum runs over the singular points of the vector field $\widetilde{X}$ on the fibre $V^F_t$.
\end{definition}

The Poincar\'e--Hopf index depends on the choice of a smoothing, but does not depend on the choice of $t$ and of the extension $\widetilde{X}$ \cite[Proposition~3.4.1]{BSS09}. In particular, one has the following proposition (see also \cite[Proposition~3.4.1]{BSS09}).

\begin{proposition} \label{prop:PH=chi}
Let the vector field $X$ be transversal to the link $K=V \cap S_\eps$ of the singularity $(V,0)$. Then
\[
\indPH^F(X;V,0) = \chi(V^F_t).
\]
\end{proposition}

Now let $(V,0)$ be the germ of a complete intersection with an isolated singularity at the origin. Then there is an essentially unique smoothing of $(V,0)$, since the base space of the semi-universal deformation is smooth (see \cite[Theorem~7.2.22]{Greuel20}). Therefore we can write in this case
$\indPH(X;V,0):=\indPH^F(X;V,0)$, where $F$ is the unique smoothing of $(V,0)$. In this case we have:

\begin{proposition} \label{prop:PH=GSV}
For a vector field $X$  on an {\icis}\ $(V,0)$  we have
\[
\indPH(X;V,0) = \indGSV(X;V,0).
\]
\end{proposition}

Therefore the Poincar\'e--Hopf index relative to a smoothing is called the GSV index relative to a smoothing in \cite{BSS09}. 

Seade defined in this way an index for a singular point of a vector field on a complex analytic surface with a smoothable normal Gorenstein singularity \cite{Seade87}.

The following corollary of Proposition~\ref{prop:PH=chi} and Proposition~\ref{prop:PH=GSV} is proved in \cite[Proposition~1.4]{SS98}.

\begin{proposition} \label{prop:muX}
For a vector field $X$  on an {\icis} \ $(V,0)$  we have
\[ 
\indGSV(X;V,0) = \indrad(X;V,0)+(-1)^n\mu ,
\]
where $\mu$ is the Milnor number of  $(V,0)$.
\end{proposition}

Similarly, for a 1-form one can prove \cite[Proposition~2.8]{EGS04}:
\begin{proposition} \label{prop:mu}
For a 1-form $\omega$  on an {\icis} \ $(V,0)$  we have
\[ 
\indGSV(\omega;V,0) = \indrad(\omega;V,0)+\mu .
\]
\end{proposition}

There is a generalisation of this index to vector fields on germs of complete intersections with non-isolated singularities \cite{BSS05c}, see also \cite[Section~3.5]{BSS09}. Let $(V,0)$ be the germ of a complete intersection defined by a map germ $F=(f_1, \ldots , f_k): (\CC^{N},0) \to (\CC^k,0)$. We assume that a neighbourhood of the origin 
in $\CC^n$ permits a Whitney stratification adapted to $V$ and satisfying the Thom condition ($a_f$) \cite[Definition~4.4.1]{Trotman20}. This holds, in particular, if $(V,0)$ is an \icis\ (say, a hypersurface, i.e.\ $k=1$). Let $X$ be a stratified vector field on $(V,0)$ with an isolated singularity at the origin. Then one can define in a similar way as above a Poincar\'e--Hopf or GSV index $\indPH(X;V,0) = \indGSV(X;V,0)$. For details see \cite[Section~3.5]{BSS09}.

\subsection{Homological index} \label{sect:Hom}
Let $(V,0) \subset (\CC^N,0)$ be a germ of a complex analytic variety of pure dimension $n$ with an isolated singular point at the origin. Let $X$ be a complex analytic vector field tangent to $(V,0)$ with an isolated singular point at the origin and let $\omega$ be a holomorphic 1-form on $(V,0)$ with an isolated singularity at the origin. We shall define an index of $X$ and $\omega$ in a homological way.

For this purpose, we consider the module $\Omega^k_{V,0}$ of germs of holomorphic $k$-forms on $(V,0)$. It is defined as follows. Let $I_{V,0} \subset \calO_{\CC^N,0}$ be the ideal of germs of holomorphic functions vanishing on $(V,0)$. Consider the $\calO_{\CC^N,0}$-module $\Omega^k_{\CC^N,0}$ of germs of holomorphic $k$-forms on $(\CC^N,0)$. Then 
\[ \Omega^k_{V,0}=\Omega^k_{\CC^N,0}/\{ f \cdot \Omega^k_{\CC^N,0} + df \wedge \Omega^{k-1}_{\CC^N,0} \, : \, f \in I_{V,0} \}.
\]
We consider two Koszul complexes:
\begin{eqnarray}
(\Omega^\bullet_{V,0}, X)& : & 
0 \longleftarrow \calO_{V,0}
\stackrel{X}{\longleftarrow} \Omega^1_{V,0}
\stackrel{X}{\longleftarrow} ... \stackrel{X}{\longleftarrow}
\Omega^n_{V,0} \longleftarrow 0\,, 
\label{eqn:hom_vector}\\
(\Omega^\bullet_{V,0}, \wedge\omega) & : &
0 \longrightarrow \calO_{V,0} \stackrel{\wedge \omega}{\longrightarrow}
\Omega^1_{V,0}
\stackrel{\wedge \omega}{\longrightarrow} ...
\stackrel{\wedge \omega}{\longrightarrow} \Omega^n_{V,0} \longrightarrow 0\,.\label{eqn:hom_1-form}
\end{eqnarray}
For the first complex $(\Omega^\bullet_{V,0}, X)$, the arrows are given by contraction with the vector field $X$. For the second complex $(\Omega^\bullet_{V,0}, \wedge\omega)$, the arrows are given by the exterior product with the 1-form $\omega$. The second complex is the dual of the first one. It was used by G.~M.~Greuel in \cite{Greuel75}. The sheaves $\Omega^i_{V,0}$ are coherent sheaves and the cohomology  sheaves of the complexes are concentrated 
at the origin and hence finite dimensional. 

\begin{remark} 
If $(V,0)=(W,0)\times(\CC,0)$ and
$X=\frac{\partial{\ }}{\partial t}$, where $t$ is the coordinate on $\CC$, the homology groups of the complex (\ref{eqn:hom_vector}) are trivial. 
This implies that, if $X$ has an isolated singular point at the origin, the homology groups $H_i(\Omega^{\bullet}_{V,0},X)$ of the complex (\ref{eqn:hom_vector}) are finite dimensional. even if $V$ has non-isolated  singularities.
\end{remark}

\begin{remark} \label{rmk:holomorphic_vf}
On a germ of a complex analytic variety with an isolated singularity, holomorphic vector fields with isolated singular points always exist \cite{BG94}. This is not the case for varieties with non-isolated singularities. For example, on the surface in $\CC^3$ given by 
\[ xy(x+y)(x+zy)=0
\]
all holomorphic vector fields vanish on the $z$-axis (cf.\ \cite[Example~13.2]{Whitney}, see also \cite{EGGD}).
\end{remark}

Let us denote by $h_j(\Omega^\bullet_{V,0},X)$ and $h_j(\Omega^\bullet_{V,0},\wedge\omega)$ the dimension (as a $\CC$-vector space) of  the $j$-th homology group of the complex $(\Omega^\bullet_{V,0}, X)$ and $(\Omega^\bullet_{V,0}, \wedge\omega)$ respectively.

\begin{definition}
(a) The {\em homological index} \index{homological index! of a vector field}
$\,{\ind}_{\rm hom}(X; V,0)$ of the vector field
$X$ on $(V, 0)$ is the Euler characteristic of the
complex $(\Omega^\bullet_{V,0}, X)$:
\begin{equation}\label{eq0}
{\ind}_{\rm hom}(X; V,0)  = \sum_{j=0}^n
(-1)^j h_j(\Omega^\bullet_{V,0},X)\,.
\end{equation}

(b) The {\em homological index} \index{homological index! of a 1-form}
$\,{\ind}_{\rm hom}(\omega; V,0)$ of the 1-form
$\omega$ on $(V, 0)$ is $(-1)^n$ times the Euler characteristic of the 
complex $(\Omega^\bullet_{V,0}, \wedge\omega)$:
\begin{equation}\label{eq1}
{\ind}_{\rm hom}(\omega; V,0) = \sum_{j=0}^n
(-1)^{n-j} h_j(\Omega^\bullet_{V,0},\wedge\omega)\,.
\end{equation}
\end{definition}

The definition of the homological index of a vector field is due to G\'omez-Mont \cite{GomezMont98}. The definition was adapted to the case of a 1-form in \cite{EGS04}. Both indices satisfy the law of conservation of number \cite[Theorem~1.2]{GomezMont98}, \cite{GG02}. 

\begin{theorem} Let $(V,0)$ be an isolated complete intersection singularity. Then 
\begin{eqnarray*}
{\ind}_{\rm hom}(X; V,0) & = & \indGSV(X;V,0),\\
{\ind}_{\rm hom}(\omega; V,0) & = & \ind_{\rm GSV} \, (\omega;V,0).
\end{eqnarray*}
\end{theorem}
For a vector field $X$, this was proved in \cite[Theorem~3.5]{GomezMont98} for a hypersurface singularity and by H.-Ch.~Graf von Bothmer, G\'omez-Mont and the first author in \cite[Theorem~2.4]{BEG08} for a complete intersection singularity. For the proof for a 1-form $\omega$ see \cite[Theorem~3.2 (iii)]{EGS04}. In \cite{Alex14, Alex16}, methods of computation of the homological index are given for Cohen--Macaulay curves, graded normal surfaces, and complete intersections.

\begin{example} Let $(C,0)$ be a curve singularity and let $(\overline{C},\overline{0})$ be its normalization. Let $\tau = \dim {\rm Ker}(\Omega^1_{c,0} \to  \Omega^1_{\overline{C},\overline{0}})$ and $\lambda=\dim \omega_{C,0}/c(\Omega^1_{C,0})$, where $\omega_{C,0}$ is the dualizing module of Grothendieck and $c: \Omega^1_{C,0} \to \omega_{C,0}$ is the class map (see \cite{BG80}). A Milnor number $\mu(f)$ of a function $f$ on a curve singularity was introduced for curves in $\CC^3$ in \cite{Goryunov00} and for the general case in \cite{MvS01}. In a similar way, one can define a Milnor number for an analytic 1-form $\omega$ with a isolated singular point on $(C,0)$, namely 
\[ \mu(\omega) := \dim \omega_{C,0}/\omega \wedge \calO_{C,0}. \]
Then one has
\[ \mu(\omega)= \indhom(\omega;C,0) + \lambda - \tau.
\]
For other formulas for $\mu(f)$ see \cite{N-BT08}.
\end{example}

It follows from Proposition~\ref{prop:mu} that, for a 1-form on an {\icis} $(V,0)$, the difference 
\[ {\ind}_{\rm hom}(\omega; V,0)-{\ind}_{\rm rad}(\omega; V,0)
\]
between the homological index and the radial index is equal to the Milnor number of the singularity. This difference is also defined for the germ of a complex analytic space of pure dimension $n$ with an isolated singular point at the origin. By \cite[Proposition~4.1]{EGS04}, it does not depend on the 1-form $\omega$. Therefore one can consider this difference as a generalized Milnor number of the singularity $(V,0)$. In \cite{EGS04}, this invariant is computed for arbitrary curve singularities and compared with the Milnor number introduced by R.~Buchweitz and G.~M.~Greuel \cite{BG80} for these singularities. See \cite{Seade17} for an interesting question about this Milnor number for normal surface singularities. See \cite{Seade19} for a survey on relations between Milnor numbers and indices of vector fields and 1-forms on singular varieties.

\subsection{Euler obstruction} \label{sect:Eu}
In this section, we define the local Euler obstruction of a singular point of a vector field or a 1-form. The idea goes back to MacPherson who defined in \cite{MacPherson74} the Euler obstruction of a singular point of a complex analytic variety. The definition of the Euler obstruction of a singular point of a vector field for the case of radial extensions is due to Brasselet and Schwartz \cite{BS81}. It seems that a definition for arbitrary vector fields is not present in the literature.

In order to introduce this notion, we need the notion of the Nash transformation of a germ of a singular variety. Let $(V,0) \subset (\CC^N,0)$ be the germ of a purely $n$-dimensional complex analytic variety. We assume that $V$ is a representative of $(V,0)$ defined in a suitable neighbourhood $U$ of the origin in $\CC^N$.  Let $G(n,N)$ be the Grassmann manifold of $n$-dimensional vector subspaces of $\CC^N$. Let $V_{\rm reg}$ be the non-singular part of $V$. There is a natural map $\sigma: V_{\rm reg} \to U \times G(n,N)$ which is defined by $\sigma(z)=(z,T_zV_{\rm reg})$. The {\em Nash transform}\index{Nash transform}\index{transform!Nash} $\widehat{V}$ is the closure of the image
${\rm Im}\,\sigma$ of the map $\sigma$ in $U \times G(n,N)$. It is a usually singular analytic variety. There is the natural base point map $\nu: \widehat{V} \to V$. Let $\widehat{V}':=\widehat{V} \setminus \nu^{-1}(V \setminus V_{\rm reg})$. Then the restriction $\nu|_{\widehat{V}'}$ maps $\widehat{V}'$ biholomorphically to $V_{\rm reg}$.

The {\em Nash bundle}\index{Nash bundle}\index{bundle!Nash}  $\widehat{T}$ over $\widehat{V}$ is the pullback of the tautological bundle on the Grassmann manifold $G(n,N)$ under the natural projection map $\widehat{V} \to G(n,N)$. It is a vector bundle of rank $n$. There is a natural lifting of the Nash transformation to a bundle
map from the Nash bundle $\widehat{T}$ to the restriction of the tangent
bundle $T\CC^N$ of $\CC^N$ to $V$. This is an isomorphism of $\widehat{T}$
and $TV_{\rm reg} \subset T\CC^N$ over the regular part $V_{\rm reg}$ of $V$.

Let $V=\bigcup_{i=1}^q V_i$ be a subanalytic Whitney stratification of $V$ and let $X$ be a stratified vector field on $V$. Let $0 \in V$ be an isolated singular point of $X$ on $V$. By \cite{BS81}, the vector field $X$ has a canonical lifting to a section $\widehat{X}$ of the Nash bundle $\widehat{T}$ over the Nash transform $\widehat{V}$ without zeros outside of $\nu^{-1}(0)$. Let $\eps$ be chosen such that the vector field $X$ is defined on $B_\eps \cap V$  and does not vanish there outside of the origin where $B_\eps$ is the ball of radius $\eps$ centered at the origin.

\begin{definition}
The {\em local Euler obstruction}\index{local Euler obstruction! of a vector field}\index{Euler obstruction \see{local Euler obstruction}} ${\rm Eu}(X; V,0)$ of the vector field $X$ on $V$ at the origin is the
obstruction to extend the non-zero section $\widehat{X}$ from the
preimage of a neighbourhood of the sphere $S_\eps=\partial B_\eps$ to
the preimage of its interior, more precisely its value $($as an element of
$H^{2n}(\nu^{-1}(V\cap B_\eps),\nu^{-1}(V\cap S_\eps))$\,$)$ on the fundamental
class of the pair $(\nu^{-1}(V\cap B_\eps), \nu^{-1}(V \cap S_\eps))$.
\end{definition} 

Now let $\omega$ be a 1-form on $U$ with an isolated singular point on $V$ at the origin. Let $\eps$ be small enough such that the 1-form $\omega$ has no singular points on $V \setminus \{ 0 \}$ inside the ball $B_\eps$. The 1-form $\omega$ gives rise to a section $\widehat{\omega}$ of the dual Nash bundle $\widehat{T}^\ast$ over the Nash transform $\widehat{V}$ without zeros outside of $\nu^{-1}(0)$. The following definition was given in \cite{EGGD}.

\begin{definition} The {\em local Euler obstruction}\index{local Euler obstruction! of a 1-form} ${\rm Eu}(\omega; V,0)$
of the 1-form $\omega$ on $V$ at the origin is the obstruction to extend
the non-zero section $\widehat{\omega}$ from the preimage of a
neighbourhood of the sphere $S_\eps= \partial B_\eps$ to the preimage
of its interior, more precisely its value (as an element of the cohomology
group $H^{2n}(\nu^{-1}(V\cap B_\eps), \nu^{-1}(V \cap S_\eps))$\,) on the
fundamental class of the pair
$(\nu^{-1}(V\cap B_\eps), \nu^{-1}(V \cap S_\eps))$.
\end{definition}

We can use the definition of the local Euler obstruction of a 1-form to define the local Euler obstruction of a germ of a complex variety which was the original definition of MacPherson in \cite{MacPherson74}.

\begin{definition}
The {\em local Euler obstruction}\index{local Euler obstruction! of a germ} ${\rm Eu}(V,0)$ of the germ $(V,0)$ is the local Euler obstruction of the radial 1-form $d|r|^2$ on it, where $r$ is the distance to the origin.
\end{definition}

The word {\em local} will be usually omitted.

One has the following Proportionality Theorem due to Brasselet and Schwartz \cite{BS81}. (For a proof see also \cite[Section~8.1.1]{BSS09}, \cite{BSS05b}.) Let $V$ be a complex analytic variety with a Whitney stratification $\{ V_i \}$.

\begin{theorem}[Proportionality Theorem for vector fields] \label{thm:propX}
Let $V_i$ be a stratum of the Whitney stratification and  $x \in V_i$, let $X_i$ be a vector field on $V_i$ with an isolated singular point at $x$, and let $X$ be a radial extension of $X_i$. Then one has
\[ {\rm Eu}(X; V,x) = {\rm Eu}(V,x) \cdot \indrad(X;V,x).
\]
\end{theorem}

Note that $\indrad(X;V,x)=\ind(X_i;V_i,x)$.

There is also a Proportionality Theorem for 1-forms due to Brasselet, Seade, and Suwa \cite{BSS07}.

\begin{theorem}[Proportionality Theorem for 1-forms] \label{thm:prop1-form}
Let $V_i$ be a stratum of the Whitney stratification and $x \in V_i$, let $\omega_i$ be a 1-form on $V_i$ with an isolated singular point at $x$, and let $\omega$ be a radial extension of $\omega_i$. Then one has
\[ {\rm Eu}(\omega; V,x) = {\rm Eu}(V,x) \cdot \indrad(\omega;V,x).
\]
\end{theorem}

It follows from Theorem~\ref{thm:prop1-form} that, if $V_i$ and $V_j$ are strata of the Whitney stratification with $V_i \subset \overline{V}_j$, then the local Euler obstruction ${\rm Eu}(\overline{V}_j,p)$ at any point $p \in V_i$ does not depend on $p$. It will be denoted by ${\rm Eu}(V_j, V_i)$. It is equal to the local Euler obstruction ${\rm Eu}(N_{ij},p)$ of a normal slice $N_{ij}$ of the variety $\overline{V}_j$ to the stratum $V_i$ at the point $p$ \cite[Section~3]{BLS00}. If $V_i \not\subset \overline{V}_j$, we assume ${\rm Eu}(V_j, V_i)$ to be equal to zero.

In \cite{BMPS04}, the notion of the local Euler obstruction of a holomorphic function $f$ with an isolated critical point on $(V,0)$ was introduced. It is defined as follows. Let $f$ be a holomorphic function defined in $U$ with a isolated singular point on $V$ at the origin. Let $\eps>0$
be small enough such that the function $f$ has no singular points on
$V\setminus\{0\}$ inside the ball $B_\eps$. Let ${\rm grad}\,f$ be
the gradient vector field of $f$ as defined in \ref{sect:GSV}.
Since $f$ has no singular
points on $V \setminus \{0\}$ inside the ball $B_\eps$, the angle of
${\rm grad}\, f(x)$ and the tangent space $T_xV_i$ to a point
$x \in V_i \setminus \{ 0 \}$ is less than $\pi/2$. Denote by
$\zeta_i(x) \neq 0$ the projection of ${\rm grad}\, f(x)$ to the tangent
space $T_xV_i$. The vector field on $V  \setminus \{ 0 \}$ which is equal to $\zeta_i$ on $V_i$ is, in general, not continuous.
It is shown in \cite{BMPS04} that the vector fields $\zeta_i$ can be glued together to obtain a stratified vector field ${\rm grad}_V f$
on $V$ such that ${\rm grad}_V f$ is homotopic to the restriction of
${\rm grad}\,f$ to $V$ and satisfies ${\rm grad}_V f(x) \neq 0$ unless $x=0$.

\begin{definition} The {\em local Euler obstruction}\index{local Euler obstruction! of a function} ${\rm Eu}(f;V,0)$ of the function $f$ is defined to be
\[ {\rm Eu}(f;V,0) := {\rm Eu}({\rm grad}_V f; V,0).
\]
\end{definition}

The Euler obstruction of a function $f$ is up to sign the Euler obstruction of the 1-form $df$. Namely, let $\omega = df$ for the germ $f$ of a holomorphic
function on $(\CC^N,0)$. Then ${\rm Eu}(df; V,0)$ differs from the
Euler obstruction ${\rm Eu}(f;V,0)$ of the function $f$ by the sign
$(-1)^n$. 
E.g., for the function $f(z_1, \ldots , z_n)=z_1^2 + \ldots + z_n^2$ on $\CC^n$ the obstruction
${\rm Eu}(f;\CC^n,0)$ is the index of the vector field
$\sum_{i=1}^n \overline{z}_i \partial/ \partial z_i$ (which is equal
to $(-1)^n$), but the obstruction ${\rm Eu}(df; \CC^n,0)$ is the index of
the (holomorphic) 1-form $\sum_{i=1}^n z_i dz_i$ which is equal to 1.

Denote by $M_f=M_{f,t_0}$ the Milnor fibre of $f$, i.e.\ the
intersection $V\cap B_\eps(0) \cap f^{-1}(t_0)$ for a (regular) value $t_0$ of $f$
close to $0$. In \cite[Theorem~3.1]{BMPS04} the following result is proved.

\begin{theorem}[Brasselet, Massey, Parameswaran, Seade] \label{theo:BMPS}
Let $f: (V,0) \to (\CC,0)$ have an isolated singularity at $0 \in V$. Then
$$
{\rm Eu}(V,0) = \left( \sum_{i=1}^q
\chi(M_f \cap V_i) \cdot {\rm Eu}(V,V_i) \right) + {\rm Eu}(f; V,0).
$$
\end{theorem}

The Euler obstruction of a vector field or 1-form can be considered as an index. In particular, it satisfies the law of conservation of
number (just as the radial index). Moreover, on a smooth variety the
Euler obstruction and the radial index coincide. This implies the
following statement (cf.\ Theorem~\ref{theo:BMPS}). We set
$\overline{\chi}(Z):=\chi(Z)-1$ and call it the {\em reduced} (modulo
a point) Euler characteristic of the topological space $Z$ (though,
strictly speaking, this name is only correct for a non-empty space $Z$).

\begin{proposition} \label{prop4}
Let $(V,0) \subset (\CC^N,0)$ have an isolated singularity at the
origin and let $\ell: \CC^N \to \CC$ be a generic linear function. Then
$$
\indrad(\omega;V,0) - {\rm Eu}(\omega; V,0) =
\indrad(d\ell;V,0) = (-1)^{n-1}\,\overline{\chi}(M_\ell),
$$
where $M_\ell$ is the Milnor fibre of the linear function $\ell$ on $V$.
In particular $$
{\rm Eu}(df; V,0) = (-1)^n (\chi(M_\ell) - \chi(M_f)).
$$
\end{proposition}

For a stratum $V_i$ of the Whitney stratification $\bigcup_{i=0}^q V_i$, $V_0=\{ 0 \}$, of $V$, let $N_i$ be the normal slice in the variety $V$ to the stratum $V_i$ at a point of the stratum $V_i$  and let 
\[ n_i=\indrad(d\ell; N_i,0)= (-1)^{\dim N_i-1}\,\overline{\chi}(M_{\ell |_{N_i}}) 
\]
be the radial index of a generic (non-vanishing) 1-form $d\ell$ on $N_i$.
In \cite[Theorem~4]{EGGD}, the following theorem was proved.

\begin{theorem} \label{thm:radEu}
One has
\[
\indrad(\omega; V,0) = \sum_{i=0}^q n_i \cdot {\rm Eu}(\omega; \overline{V}_i,0).
\]
\end{theorem}

The strata $V_i$ of $V$ are partially ordered: $V_i \prec V_j$
(we shall write $i \prec j$) iff $V_i \subset \overline{V_j}$ and
$V_i \neq V_j$; $i \preceq j$ iff $i \prec j$ or $i=j$.
In \cite[Corollary~1]{EGGD}, an ``inverse''  of the formula of Theorem~\ref{thm:radEu} was written in the case when the variety $V$ is irreducible and $V=\overline{V}_q$. 
Let $n_{ij}$ ($i \preceq j$) be the index of a generic 1-form $d\ell$ on the normal slice $N_{ij}$: $n_{ij}=  (-1)^{\dim N_{ij}-1}\,\overline{\chi}(M_{\ell |_{N_{ij}}})$ (in particular $n_{ii}=1$) and let  $m_{ij}$ be the (M\"obius) inverse of the function $n_{ij}$ on the
partially ordered set of strata, i.e.
$$
\sum_{i \preceq j \preceq k} n_{ij}m_{jk} = \delta_{ik}.
$$

\begin{corollary} \label{cor:Eurad}
One has
$$
{\rm Eu}_{V,0}\, \omega =
\sum_{i=0}^q m_{iq} \cdot {\ind}_{\rm rad} \, (\omega;\overline{V_i},0).
$$
\end{corollary}

In \cite{DGJoS}, another proof of Corollary~\ref{cor:Eurad} is given and Theorem~\ref{thm:radEu} and this corollary are applied to give an alternative proof of Theorem~\ref{theo:BMPS}.

Let $V$ be an affine variety. In \cite{STV05a}, a {\em global Euler obstruction}\index{global Euler obstruction} was defined for the variety $V$ as the obstruction to extend a radial vector field, defined outside of a sufficiently large compact subset, to a non-zero section of the Nash bundle.

In \cite{Mas20}, the local Euler obstruction was investigated in terms of constructible sheaves and characteristic cycles.

\subsection{Algebraic, analytic, and topological formulas} \label{subsec:algebraic_form}
In Section~\ref{subsec:algebraic_form-smooth}, we discussed algebraic formulas for the index of an analytic vector field or an analytic 1-form on a smooth manifold. It is natural to try to look for analogues of such formulas for vector fields and 1-forms on singular varieties. The homological index opens the way for such formulas. 

In \cite[Theorem~1]{GomezMont98}, G\'omez-Mont proved an algebraic formula for the homological index of a vector field with an isolated singularity in the ambient space on an isolated hypersurface singularity $(V,0)$. O.~Klehn \cite{Klehn03} generalized this formula to the case that the vector field has an isolated singularity on the hypersurface singularity, but not necessarily in the ambient space. Graf von Bothmer, G\'omez-Mont  and the first author  \cite{BEG08} gave formulas
to compute the homological index in the case when  $(V,0)$ is an isolated complete intersection singularity. L.~Giraldo, G\'omez-Mont, and P.~Marde\v{s}i\'c \cite{GGM02} studied the homological index of vector fields tangent to hypersurfaces with non-isolated singularities.

G\'omez-Mont and P.~Marde\v{s}i\'c derived algebraic formulas for the index of a real vector field with an algebraically isolated singular point at the origin tangent to a real analytic hypersurface with an algebraically isolated singularity at the origin as well \cite{GM97, GM99, Mar14}. The index is expressed as the signature of a certain non-degenerate quadratic form for an even-dimensional hypersurface and as the difference between the signatures of two such forms in the odd-dimensional case. 

O.~Klehn \cite{Klehn05} proved that the GSV index of a holomorphic vector field $X$ on an {\icis} $(V,0)$ coincides with the dimension of a certain explicitly constructed vector space, if $X$ is deformable in a certain sense and $V$ is a curve. Moreover, he gave a signature formula for the real GSV index in the corresponding real analytic case generalizing the Eisenbud--Levine--Khimshiashvili formula.

Now let $\omega$ be the restriction of a holomorphic 1-form $\omega = \sum_{i=1}^{N}A_i(z)dz_i$ to an {\icis} $(V,0)$ given by a mapping $f=(f_1, \ldots, f_{N-n}): (\CC^{N},0)\to(\CC^{N-n},0)$. Assume that $\omega$ has an isolated singular point at the origin (on $(V,0)$). Let $I$ be the ideal
generated by $f_1$, ..., $f_{N-n}$ and the $(N-n+1) \times (N-n+1)$-minors of
the matrix
$$\begin{pmatrix} \frac{\partial f_1}{\partial z_1} & {\cdots } &
   \frac{\partial f_1}{\partial z_{N}} \\ {\vdots} & {\cdots} &
{\vdots} \\ \frac{\partial
f_{N-n}}{\partial z_1} & \cdots &
   \frac{\partial f_{N-n}}{\partial z_{N}} \\
A_1 & \cdots & A_{N}
\end{pmatrix}\,.
$$
In \cite{EG01, EG03a}, the authors proved the following formula.
\begin{theorem}\label{thm:formulaGSV}
One has
$$
{\rm ind}_{\rm GSV} \, (\omega;V,0) = \dim {\mathcal O}_{\CC^{N},0}/I.
$$
\end{theorem}
It generalizes the L\^e--Greuel formula \cite{Greuel75, Le74} for the differential of a function. (Note that there is a minor mistake in the proof of this theorem in \cite{EG03a} which is corrected in \cite{EGBLMS}.) In \cite{EGMZ}, the authors constructed quadratic forms on the algebra ${\mathcal O}_{\CC^{N},0}/I$ and on the space $\Omega^n_{V,0}/\omega \wedge \Omega^{n-1}_{V,0}$ generalizing the Eisenbud--Levine--Khimshiashvili quadratic form defined for smooth $V$.

In~\cite{Esterov05b, Esterov05, Esterov06} A.Esterov gave formulas for the index of a 1-form on an {\icis} in terms of Newton diagrams of the components under certain genericity conditions.

T.~Gaffney \cite{Gaffney05} described connections between the GSV index of $\omega$ and the multiplicity of pairs of certain modules. 

There are several generalizations of the residue formula of Section~\ref{subsec:algebraic_form-smooth}. P.~F.~Baum and R.~Bott \cite{BB70} considered residues of meromorphic vector fields on compact complex manifolds. An integral formula for the GSV index of a holomorphic vector field on an {\icis}  was given by D.~Lehmann, M.~Soares, and Suwa in \cite{LSS95}. For generalizations and related results see \cite{LS95, Suwa95, Suwa98, Suwa02, Suwa14}.  Klehn generalized the residue formula for the GSV index to holomorphic 1-forms on an isolated surface singularity \cite{Klehn02}. T.~Honda and Suwa \cite{HS98} studied residue formulas for meromorphic functions on surfaces.
The authors \cite{EG04} considered indices of meromorphic 1-forms on complete intersections with isolated singularities. 

In \cite{Gusein-Zade84}, a topological formula for the index of a gradient vector field ${\rm grad}\, g$ of an analytic function $g :(\CC^n,0) \to (\CC,0)$ taking real values on $\RR^n \subset \CC^n$ is given. It expresses the index of the gradient vector field on $\RR^n$ in terms of signatures of certain quadratic forms on the middle homology groups of specific Milnor fibres of the germ $g$. This formula was conjectured in \cite{Arnold78} and also proved in \cite{Varchenko85}. In \cite{EG99}, a generalization of such a formula for the radial index of a gradient vector field on an algebraically isolated real analytic {\icis} in $\CC^{N}$ was obtained.

In \cite{Seade95, SS96} formulas are given evaluating the GSV index of a singular point of a vector field on an isolated hypersurface or complete intersection singularity in terms of a resolution of the singularity.

\subsection{Determinantal singularities} \label{sect:det}
The GSV index is only defined for ICIS. In the case of ICIS, the indices introduced above are best understood. In this subsection, we consider the next more general class, namely the class of determinantal singularities. An approach studying indices of 1-forms on such singularities was started in \cite{EGSteklov}. We give the basic definitions and facts following this paper.

Let $\Mmn \cong \CC^{mn}$ be the space of $m \times n$-matrices with complex entries. 

\begin{definition}
Let $t$ be an integer with $1 \leq t \leq \min(m,n)$. The {\em generic determinantal variety of type $(m,n,t)$}\index{generic determinantal variety}\index{determinantal variety!generic} is the subset 
\[ \Mmn^t := \{ A \in \Mmn \, | \, \rk(A) < t \}
\]
consisting of matrices of rank less than $t$, i.e.\ of matrices of which all ($t \times t$)-minors vanish.
\end{definition}

The variety $\Mmn^t$ has codimension $(m-t+1)(n-t+1)$ in $M_{m,n}$. It is singular. The singular locus of $\Mmn^t$ coincides with $\Mmn^{t-1}$. The singular locus of the latter one coincides with $\Mmn^{t-2}$, etc.\ (see, e.g., \cite{ACGH}). 
The representation of the variety $\Mmn^{t}$ as the union of $\Mmn^{i}\setminus\Mmn^{i-1}$, $i=1, \ldots, t$, is a Whitney stratification of $\Mmn^{t}$.

Let $U \subset \CC^N$ be an open domain and $F: U \to \Mmn$ be a holomorphic map sending $z$ to the matrix $F(z)=(f_{ij}(z))$ whose entries $f_{ij}(z)$ are complex analytic functions on $U$.

\begin{definition} A {\em determinantal variety of type $(m,n,t)$}\index{determinantal variety} is the preimage $V=F^{-1}(\Mmn^t)$ of the variety $\Mmn^t$ subject to the condition that ${\rm codim}\, V = {\rm codim}\, \Mmn^t=(m-t+1)(n-t+1)$. 
\end{definition}

The image of a generic map $F: U \to \Mmn$ may intersect the varieties $\Mmn^{i}$ for $i < t$. Therefore, it may not be avoided that $F^{-1}(\Mmn^t)$ has singularities. However, a generic map $F$ intersects the strata $\Mmn^{i} \setminus \Mmn^{i-1}$ of the variety $\Mmn^t$ transversally. This means that, at the corresponding points, the determinantal variety has ``standard'' singularities whose analytic type only depends on $i= {\rm rk}\, F(z) +1$. This inspired the following definitions of \cite[p.~114]{EGSteklov}.

\begin{definition} A point $x \in X=F^{-1}(\Mmn^t)$ is called {\it essentially non-singular}\index{essentially non-singular point}\index{point!essentially non-singular} if, at the point $x$, the map $F$ is transversal to the corresponding stratum of the variety $\Mmn^t$ (i.e., to $\Mmn^{i} \setminus \Mmn^{i-1}$ where $i= {\rm rk} \, F(x)+1$). 
\end{definition}

\begin{definition}
A germ $(V,0)\subset(\CC^N,0)$ of a determinantal variety of type $(m,n,t)$ has an {\em isolated essentially singular point} at the origin (or is an {\it essentially isolated determinantal singularity}: {EIDS})\index{essentially isolated determinantal singularity}\index{EIDS} if it has only essentially non-singular points in a punctured neighbourhood of the origin in $V$.
\end{definition}

\begin{example}
An ICIS is an example of an {EIDS}: it is an {EIDS} of type $(1,n,1)$. 
\end{example}

An essentially isolated determinantal singularity 
$(V,0) \subset (\CC^N,0)$ of type $(m,n,t)$ (defined by a map $F: (\CC^N,0) \to (\Mmn, 0)$) has an isolated singularity at the origin if and only if $N \leq (m-t+2)(n-t+2)$.

We shall consider deformations (in particular, smoothings) of an {EIDS} given by deformations of the matrix which defines the {EIDS}. Hence they are themselves determinantal ones.

Let $(V,0) \subset (\CC^N,0)$ be an EIDS defined by a map $F: (\CC^N,0) \to (\Mmn,0)$ ($V=F^{-1}(\Mmn^t)$, $F$ is transversal to $\Mmn^{i} \setminus \Mmn^{i-1}$ at all points $x$ from a punctured neighbourhood of the origin in $\CC^N$ and for all $i\le t$).

\begin{definition} An {\em essential smoothing}\index{essential smoothing}\index{smooothing!essential} $\widetilde{V}$ of the EIDS $(V,0)$ is a subvariety of a neighbourhood $U$ of the origin in $\CC^N$ defined by a perturbation $\widetilde{F}: U \to \Mmn$ of the germ $F$ transversal to all the strata $\Mmn^{i} \setminus \Mmn^{i-1}$ with $i \leq t$.
\end{definition}

A generic deformation $\widetilde{F}$ of the map $F$ defines an essential smoothing of the EIDS $(V,0)$ (according to Thom's Transversality Theorem). An essential smoothing is in general not smooth (for $N \geq (m-t+2)(n-t+2)$). Its singular locus is $\widetilde{F}^{-1}(\Mmn^{t-1})$, the singular locus of the latter one is $\widetilde{F}^{-1}(\Mmn^{t-2})$, etc. The representation of $\widetilde{V}$ as the union 
$$\widetilde{V} = \bigcup_{1 \leq i \leq t}  \widetilde{F}^{-1}(\Mmn^{i} \setminus \Mmn^{i-1})$$
is a Whitney stratification of it. An essential smoothing of an EIDS $(V,0)$ of type $(m,n,t)$ is a genuine smoothing if and only if $N < (m-t+2)(n-t+2)$.

There are three distinguished types of resolutions of the variety $\Mmn^t$.

The first one is constructed by considering $m \times n$-matrices as linear maps  $\CC^n \to \CC^m$. Let
\[
Y_1:= \{ (A,W) \in \Mmn \times G(n-t+1,n) \, | \, A(W)=0 \}.
\]
The variety $Y_1$ is smooth and connected. Its projection to the first factor defines a resolution $\pi_1: Y_1 \to \Mmn^t$ of the variety $\Mmn^t$.

Let us consider $m \times n$-matrices as linear maps $\CC^m \to \CC^n$ and let
\[
Y_2:= \{ (A,W) \in \Mmn \times G(m-t+1,m) \, | \, A^T(W)=0 \}.
\]
Then one gets a resolution $\pi_2: Y_2 \to \Mmn^t$ of the variety $\Mmn^t$.

The third natural modification is given by the Nash transform $Y_3:=\widehat{M}_{m,n}^t$. One can show that $\pi_3: Y_3 \to \Mmn^t$ is in fact a resolution of the variety $\Mmn^t$, see \cite{EGSteklov}.

Let $(V,0)=F^{-1}(\Mmn^t) \subset (\CC^N,0)$ be an EIDS and let $\omega$ be a germ of a (complex) 1-form on $(\CC^N,0)$ whose restriction to $(V,0)$ has an isolated singular point (zero) at the origin. This means that the restrictions of the 1-form $\omega$ to the strata $V_i\setminus V_{i-1}$,  $V_i:=F^{-1}(\Mmn^{i})$,
$i\le t$,
have no zeros in a punctured neighbourhood of the origin.

An essential smoothing $\widetilde{V} \subset U$ of the EIDS $(V,0)$ (in a neighbourhood $U$ of the origin in $\CC^N$) is in general not smooth. To define an analogue of the PH-index one has to construct a substitute of the tangent bundle to $\widetilde{V}$. It is possible to use one of the following two natural ways.

One possibility is to use a resolution of the variety $\widetilde{V}$ connected with one of the three resolutions of the variety $\Mmn^t$ described above. Let $\pi_k: Y_k \to \Mmn^t$ be one of the described resolutions of the determinantal variety $\Mmn^t$ and let $\overline{V}_k = Y_k \times_{\Mmn^t} \widetilde{V}$, $k=1,2,3$, be the fibre product of the spaces $Y_k$ and $\widetilde{V}$ over the variety $\Mmn^t$:
$$\diagram
 & \overline{V}_k \ar[dl]_{\Pi_k} \ar[dd] \ar[dr]  &\\
\widetilde{V} \ar[dr]_{\widetilde{F}|_{\widetilde{V}} } & & Y_k \ar[dl]^{\pi_k} \\
& \Mmn^t &
\enddiagram
$$
The map $\Pi_k : \overline{V}_k \to \widetilde{V}$ is a resolution of the variety $ \widetilde{V}$. For $k=1,2$ it is also called the {\em Tjurina transform} after \cite{Tjurina68}, see \cite{FK18}. The lifting $\omega_k:= (j \circ \Pi_k)^\ast \omega$ ($j$ is the inclusion map $\widetilde{V} \hookrightarrow U \subset \CC^N$) of the 1-form $\omega$ is a 1-form on a (non-singular) complex analytic manifold $\overline{V}_k$ without zeros outside of the preimage of a small neighbourhood of the origin. In general, the 1-form $\omega_k$ has non-isolated zeros.

\begin{definition}
The {\em Poincar\'e--Hopf index}\index{Poincar\'e--Hopf index}\index{index!Poincar\'e--Hopf} ({\em PH-index})
$\ind_{\rm PH}^k \, (\omega;V,0)$, $k=1,2,3$,
of the 1-form $\omega$ on the EIDS $(V,0) \subset (\CC^N,0)$ is the sum of the indices of the zeros of a generic perturbation $\widetilde{\omega}_k$ of the 1-form $\omega_k$ on the manifold $\overline{V}_k$ (in the preimage of a neighbourhood of the origin in $\CC^N$). 
\end{definition}

There is also the local Euler obstruction of the 1-form $\omega$. If one uses the Nash transform of the essential smoothing $\widetilde{V}$ of the EIDS $(V,0)$ instead of $V$ itself, it is called Poincar\'e--Hopf--Nash index.

\begin{definition}
The {\em Poincar\'e--Hopf--Nash index}\index{Poincar\'e--Hopf--Nash index}\index{index!Poincar\'e--Hopf--Nash} ({\em PHN-index}) $\indPHN (\omega;V,0)$ of the 1-form $\omega$ on the EIDS $(V,0)$ is the obstruction to extend the non-zero section $\widehat{\omega}$ of the dual Nash bundle $\widehat{T}^\ast$ from the preimage of the boundary $S_\eps=\partial B_\eps$ of the ball $B_\eps$ to the preimage of its interior, i.e.\ to the manifold $\overline {V}_3$, more precisely, its value (as an element of
 $H^{2d}(\Pi_3^{-1}(\widetilde{V} \cap B_\eps), \Pi_3^{-1}(\widetilde{V} \cap S_\eps))$) on the fundamental class of the pair $(\Pi_3^{-1}(\widetilde{V} \cap B_\eps), \Pi_3^{-1}(\widetilde{V} \cap S_\eps))$.
\end{definition}

In \cite[Proposition~2]{EGSteklov}, a formula relating the index $\ind_{\rm PH}^k \, (\omega;V,0)$, $k=1,2,3$ with the radial indices $\indrad(\omega; V_i,0)$, $i=0, \ldots , t$, is given. There is the following version of Theorem~\ref{thm:radEu} for the PHN-index on a determinantal variety \cite[Proposition~4]{EGSteklov}. Let $\ell: \Mmn \to \CC$ be a generic linear form and let, for $i \leq j$,
\[
n_{ij} := \indrad (d\ell; M_{m-i+1,n-i+1}^{j-i+1},0).
\]
By \cite[Proposition~3]{EGSteklov}, we have
\begin{eqnarray*}
\lefteqn{\indrad (d\ell; M_{m-i+1,n-i+1}^{j-i+1},0)} \\
 & = & (-1)^{d_{ij}-1} \overline{\chi}(M_{m-i+1,n-i+1}^{j-i+1}\cap \ell^{-1}(1)) = (-1)^{(m+n)(j-i)} \binom{m-i}{m-j},
\end{eqnarray*}
where $d_{ij}$ is the dimension of $M_{m-i+1,n-i+1}^{j-i+1}$ equal to $(m-i+1)(n-i+1)-(m-j+1)(n-j+1)$.

\begin{theorem} \label{Propnit}
One has
$$\indrad (\omega;V,0) = \sum_{i=1}^t n_{it} \indPHN (\omega; V_i,0) + (-1)^{\dim V -1} \overline{\chi}(\widetilde{V},0).$$
\end{theorem}

One can see that, for $i \leq j \leq t$, the integers $m_{ij}$ from Corollary~\ref{cor:Eurad} are given by
$$m_{ij} = (-1)^{(m+n+1)(j-i)} \binom{m-i}{m-j}.$$
The analogue of Corollary~\ref{cor:Eurad}, the inverse to Theorem~\ref{Propnit}, is the following statement, see \cite[Proposition~5]{EGSteklov}.

\begin{corollary} \label{PropPHN}
One has
$$\indPHN (\omega;V,0) = \sum_{i=1}^t m_{it} \left( \indrad (\omega; V_i,0) + (-1)^{\dim V_i} \overline{\chi}(\widetilde{V}_i,0) \right).$$
\end{corollary}

T.~Gaffney, N.~G.~Grulha, and  M.~A.~S.~Ruas \cite{GGR} proved a generalization of  Theorem~\ref{Propnit} and its inverse Corollary~\ref{PropPHN} relating it with Gaffney's multiplicities of pairs of modules.

For isolated determinantal singularities, the relations between the PH-, the PHN- and the radial indices simplify. For isolated smoothable singularities (i.e.\ for $N < (m-t+2)(n-t+2)$)  all Poincar\'e--Hopf indices (including the Poincar\'e--Hopf--Nash index) coincide and they are equal to
$$\indPH(\omega; V,0) =  \indrad (\omega; V,0) + (-1)^{\dim V} \overline{\chi}(\widetilde{V},0).$$

The paper \cite{EGSteklov} contains an algebraic formula for this index: Proposition~8 therein. To a regret, its proof is wrong.

N.~C.~Chachapoyas Siesqu\'en \cite{ChS18} studies the Euler obstruction of an EIDS and gives some formulas to calculate it. In \cite{N-BOT13} (see also \cite{N-BOT18}), a formula for the Euler obstruction of a smoothable IDS is given. The papers \cite{DGP} and \cite{ANOT}  contain results on the Euler obstruction of a function on a determinantal variety.

In \cite{RuasPe14}, codimension two determinantal varieties with isolated singularities are studied. The Milnor number is defined to be the middle Betti number of a generic fibre of the unique smoothing of such a singularity. For surfaces in $\CC^4$, a L\^e--Greuel formula for the Milnor number of the surface is proved. The Milnor number is also related to the Poincar\'e--Hopf index of the 1-form given by the differential of a generic linear projection defined on the surface. For other generalizations of the L\^e--Greuel formula  see \cite{DGAdv, C-BMSS}. For other results on Milnor numbers of essentially isolated determinantal singularities see \cite{FKZ15, BCR17}.

\section{Indices of collections of vector fields and 1-forms}\label{sec:coll}
\subsection{GSV index} \label{sect:GSVcoll}
Let $(V,0) \subset (\CC^N,0)$ be an {\icis} defined by a holomorphic map germ
$f=(f_1, \ldots , f_{N-n}): (\CC^N,0) \to (\CC^{N-n},0)$. 
Let $\{ X_j^{(i)} \}$ be a collection of vector
fields on a neighbourhood of the origin in $(\CC^N,0)$ ($i=1, \ldots , s$;
$j=1, \ldots , n-k_i+1$; $\sum k_i = n$) which are tangent to the
{\icis} $(V,0)=\{f_1= \cdots = f_{N-n}=0\} \subset (\CC^N,0)$ at non-singular
points of $V$.
We say that a point
$p\in V\setminus \{0\}$ is non-singular for the collection $\{ X_j^{(i)} \}$
 on $V$ if at least for some $i$ the vectors $X_1^{(i)}(p), \ldots , X_{n-k_i+1}^{(i)}(p)$ are
linearly independent. Suppose that the collection $\{ X^{(i)}_j\}$ has
no singular points on $V$ outside of the origin in a neighbourhood of it.
Let $U$ be a neighbourhood of the origin in $\CC^N$ where all the functions
$f_r$ ($r=1, \ldots , N-n$) and the vector fields $X_j^{(i)}$ are defined
and such that the collection $\{ X_j^{(i)} \}$ has no singular
points on $(V \cap U) \setminus \{ 0\}$. Let
$S_\delta \subset U$ be a sufficiently small sphere around the origin
which intersects $V$ transversally and denote by $K=V \cap S_\delta$
the link of the {\icis} $(V,0)$. The manifold $K$ has a natural orientation
as the boundary of a complex analytic manifold. Let $\Psi_V$ be the mapping
from $V \cap U$ to $M_{N,{\bf k}}$ (for the definition of $M_{N,{\bf k}}$ see Section~\ref{sect:coll}) which sends a point
$x \in V \cap U$ to the collection of $N \times (N-k_i+1)$-matrices
$$
\{ ( \mbox{grad}\, f_1(x), \ldots , \mbox{grad}\, f_{N-n}(x),
X_1^{(i)}(x), \ldots , X_{n-k_i+1}^{(i)}(x))\}, \quad i=1,\ldots, s.
$$
Here
$\mbox{grad}\, f_r$ is the gradient vector field of $f_r$ defined in
Section~\ref{sect:GSV}. Its restriction $\psi_V$ to the link $K$ maps $K$ to the subset
$W_{N,{\bf k}}$.

\begin{definition}
The
{\em GSV index} \index{GSV index! of a collection of vector fields}
$\indGSV(\{X^{(i)}_j \}; V,0)$ of the
collection of vector fields $\{ X^{(i)}_j \}$ on the {\icis} $(V,0)$
is the degree of the mapping $\psi_V : K \to W_{N, {\bf k}}$, or, equivalently, the intersection number of the germ of the image of the map $\Psi_V$ with the variety $D_{N,{\bf k}}$.
\end{definition}

For $s=1$, $k_1=n$, this index is the GSV index of a vector field on an {\icis} defined in Section~\ref{sect:GSV}.

Let $V \subset \CC\PP^N$ be an $n$-dimensional complete intersection with isolated singular points which is defined by homogeneous polynomials $f_1, \ldots , f_{N-n}$ in $(N+1)$ variables. Let $\{ X^{(i)}_j \}$ be
a collection of continuous vector fields on $\CC\PP^N$ which are tangent
to $V$. Let $\widetilde V$ be a smoothing of the complete intersection $V$,
i.e.\ $\widetilde V$ is defined by $N-n$ homogeneous polynomials $\widetilde f_1, \ldots, 
\widetilde f_{N-n}$ which are small perturbations of the functions $f_i$ and $\widetilde V$
is smooth. As in Section~\ref{sect:PH}, one can define approximations $\{ \widetilde{X}^{(i)}_j \}$ of the vector fields $\{ X^{(i)}_j \}$ which are tangent to $\widetilde{V}$. Then one has the following analogue of Theorem~\ref{thm:Chernsec}.

\begin{theorem} \label{thm:ChernX}
One has
$$
\sum_{p \in V} \indGSV(\{X^{(i)}_j \}; V,p) = \langle \prod_{i=1}^s
c_{k_i}(T \widetilde{V}), [ \widetilde{V} ] \rangle,
$$
where $\widetilde{V}$ is a smoothing of the complete intersection $V$.
\end{theorem}

Now let
$\{ \omega^{(i)}_j\}$ be a collection of (continuous) 1-forms on a
neighbourhood of the origin in $(\CC^N,0)$ with $i=1, \ldots, s$,
$j=1, \ldots , n-k_i+1$, $\sum k_i = n$. We say that a point
$p\in V\setminus \{0\}$ is non-singular for the collection
$\{\omega^{(i)}_j\}$ on $V$ if at least for some $i$ the restrictions
of the 1-forms $\omega^{(i)}_j(p)$, $j=1, \ldots , n-k_i+1$, to the tangent space $T_p V$ are
linearly independent. Assume that the collection $\{ \omega^{(i)}_j\}$ has no singular points on $V$ in a punctured neighbourhood of the origin. As above,
let $U$ be a neighbourhood of the origin in $\CC^N$ where all the functions
$f_r$ ($r=1, \ldots , N-n$) and the 1-forms $\omega_j^{(i)}$ are defined
and such that the collection $\{ \omega_j^{(i)} \}$ has no singular
points on $(V \cap U) \setminus \{ 0\}$. Let $S_\delta \subset U$ be a sufficiently small sphere around the origin
As above, let $K=V\cap S_\delta$ be the link of the {\icis} $(V,0)$.
Let $\Psi_V$ be the mapping
from $V \cap U$ to $M_{n,{\bf k}}$ which sends a point
$x \in V \cap U$ to the collection of $N \times (N-k_i+1)$-matrices
$$
\{ (df_1(x), \ldots , df_{N-n}(x), \omega_1^{(i)}(x), \ldots ,
\omega_{n -k_i+1}^{(i)}(x)) \}, \quad i=1, \ldots, s.
$$
Its restriction $\psi_V$ to the link $K$ maps $K$ to the subset
$W_{N,{\bf k}}$.

\begin{definition}
The {\em GSV index} \index{GSV index! of a collection of 1-forms} $\indGSV(\{\omega^{(i)}_j\};V,0)$ of the
collection of 1-forms $\{ \omega^{(i)}_j \}$ on the {\icis} $(V,0)$
is the degree of the mapping $\psi_V : K \to W_{N, {\bf k}}$, or, equivalently, the intersection number of the germ of the image of the
mapping $\Psi_V$ with the variety $D_{N, {\bf k}}$.
\end{definition}

For $s=1$, $k_1=n$, this index is the GSV index of a 1-form on an {\icis} defined in Section~\ref{sect:GSV}.

Let $V \subset \CC\PP^N$ be an $n$-dimensional complete intersection with isolated singular points which is defined by homogeneous polynomials $f_1, \ldots , f_{N-n}$ in $(N+1)$ variables. Let $L$ be a complex line
bundle on $V$ and let $\{ \omega^{(i)}_j \}$ be a collection of continuous
1-forms on $V$ with values in $L$. This means that the forms
$\omega_j^{(i)}$ are continuous sections of the vector bundle
$T^\ast V \otimes L$ outside of the singular points of $V$.
Since, in
a neighbourhood of each point $p$, the vector bundle $L$ is trivial,
one can define the index $\indGSV(\{\omega^{(i)}_j\};V,p)$ of the collection of 1-forms $\{ \omega^{(i)}_j \}$ at the point $p$  as above.
Let $\widetilde V$ be a smoothing of the complete intersection $V$.
By using, e.g., the pull back along a projection of $\widetilde{V}$ to $V$, one can consider $L$ as a line bundle
on $\widetilde{V}$ as well. The collection $\{\omega^{(i)}_j \}$ of 1-forms can also be extended
to a neighbourhood of $V$ in such a way that it will define a collection
of 1-forms on the smoothing $\widetilde{V}$ (also denoted by
$\{\omega^{(i)}_j \}$) with isolated singular points.
The sum of the indices of the collection $\{\omega^{(i)}_j \}$
on the smoothing $\widetilde{V}$ of $V$ in a neighbourhood of the point
$p$ is equal to the index $\indGSV(\{\omega^{(i)}_j\};V,p)$.

One has the following analogue of Theorem~\ref{thm:ChernX} for
1-forms.

\begin{theorem} \label{sec8theo1}
One has
$$
\sum_{p \in V} \indGSV(\{\omega^{(i)}_j\};V,p) = \langle \prod_{i=1}^s
c_{k_i}(T^\ast\widetilde{V} \otimes L), [ \widetilde{V} ] \rangle,
$$
where $\widetilde{V}$ is a smoothing of the complete intersection $V$.
\end{theorem}

Now let $(V,0)$ be an {\icis} defined by a holomorphic map germ
$f=(f_1, \ldots , f_{N-n}): (\CC^N,0) \to (\CC^{N-n},0)$ as above.
Let $\{ \omega_j^{(i)}\}$
($i=1, \ldots, s$; $j=1, \ldots , n-k_i+1$) be a collection of
1-forms on a neighbourhood of the origin in $\CC^N$ without singular
points on $V\setminus\{0\}$ in a neighbourhood of the origin. We now assume that all the 1-forms $\omega_j^{(i)}$ are complex analytic.

Let $I_{V,\{ \omega_j^{(i)}\}}$ be
the ideal in the ring ${\mathcal O}_{\CC^N,0}$ generated by the functions
$f_1, \ldots , f_{N-n}$ and by the $(N-k_i+1) \times (N-k_i+1)$ minors
of all the matrices
$$
(df_1(x), \ldots , df_{N-n}(x), \omega_1^{(i)}(x), \ldots,
\omega_{n-k_i+1}^{(i)}(x))
$$
for all $i=1, \ldots, s$. Then we have the following algebraic formula similar to that of Theorem~\ref{thm:formulaGSV} (see \cite{EGBLMS}).

\begin{theorem} \label{sec8theo2}
$$
\indGSV(\{\omega^{(i)}_j\};V,0) =
\dim_\CC {\mathcal O}_{\CC^N,0}/I_{V,\{\omega_j^{(i)}\}}.
$$
\end{theorem}

\begin{remark}
In the case of collections of vector fields, the map $\Psi_V$ is not complex analytic, whereas it is complex analytic in the case of collections of 1-forms. This is the reason that 
a formula similar to that of
Theorem~\ref{sec8theo2} does not exist for collections of vector fields.  Moreover,
in some cases this index can be negative (see e.g.
\cite[Proposition~2.2]{GSV91}).
\end{remark}

\subsection{Chern obstruction} \label{sect:Chern}
We now consider a generalization of the notion of the Euler obstruction to
collections of 1-forms corresponding to different Chern numbers.

Let $(V,0)\subset(\CC^N,0)$ be the germ of a purely $n$-dimensional
reduced complex analytic variety at the origin. It can have 
a non-isolated singularity at the origin. Let ${\bf k}=\{k_i\}$, $i=1,\ldots, s$, be
a fixed partition of $n$ (i.e., $k_i$ are positive integers,
$\sum\limits_{i=1}^s k_i=n$). Let $\{\omega^{(i)}_j\}$ ($i=1,\ldots, s$,
$j=1,\ldots, n-k_i+1$) be a collection of germs of 1-forms on $(\CC^N, 0)$
(not necessarily complex analytic; it is sufficient  that they are continuous). Let $\eps>0$ be small enough so that there is
a representative $V$ of the germ $(V,0)$ and representatives
$\omega^{(i)}_j$ of the germs of 1-forms inside the ball
$B_\eps \subset\CC^N$.

\begin{definition}
A point $p\in V$ is called a {\em special} \index{special point} point of the collection
$\{\omega^{(i)}_j\}$ of 1-forms on the variety $V$ if there exists
a sequence $\{p_m\}$ of points from the non-singular part $V_{\rm reg}$
of the variety $V$ such that the sequence $T_{p_m}V_{\rm reg}$ of the
tangent spaces at the points $p_m$ has a limit $L$ $($in $G(n,N)$ $)$ as $m$ tends to infinity and the
restrictions of the 1-forms $\omega^{(i)}_1$, \dots, $\omega^{(i)}_{n-k_i+1}$
to the subspace $L\subset T_p\CC^N$ are linearly dependent for each
$i=1, \ldots, s$. The collection $\{\omega^{(i)}_j\}$ of 1-forms has an
{\em isolated special point} \index{isolated special point} on $(V,0)$ if it has no special points on
$V$ in a punctured neighbourhood of the origin.
\end{definition}

For the case $s=1$ (and therefore $k_1=n$), i.e. for one 1-form
$\omega$, we discussed the notion of a {\em singular} point of the 1-form
$\omega$ on $V$ in Section~\ref{sect:rad}. One can easily see that a special point of the 1-form
$\omega$ on $V$ is singular, but not vice versa. (E.g.\ the origin is a
singular point of the 1-form $dx$ on the cone $\{x^2+y^2+z^2=0\}$, but
not a special one.) On a smooth variety these two notions coincide.

Let
$$
{\mathcal L}^{\bf k}= \prod\limits_{i=1}^s \prod\limits_{j=1}^{n-k_i+1}
(\CC^N_{ij})^\ast
$$
be the space of collections of linear functions on $\CC^N$ (i.e.\ of
1-forms with constant coefficients). Then one can show \cite[Proposition~1.1]{EGBrasselet}
that there exists an open and dense set $U \subset {\mathcal L}_{\bf k}$ such that
each collection $\{\ell^{(i)}_j\} \in U$ has only isolated special points
on $V$ and, moreover, all these points belong to the smooth part
$V_{\rm reg}$ of the variety $V$ and are non-degenerate (see Section~\ref{sect:coll} for the notion of a non-degenerate singular point).
This implies the following proposition (see \cite[Corollary~1.1]{EGBrasselet}).

\begin{proposition} \label{prop:Morse}
Let $\{\omega^{(i)}_j\}$ be a collection of 1-forms on $V$ with an
isolated special point at the origin. Then there exists a deformation
$\{\widetilde \omega^{(i)}_j\}$ of the collection $\{\omega^{(i)}_j\}$
whose special points lie in $V_{\rm reg}$ and are non-degenerate. Moreover,
as such a deformation one can use $\{\omega^{(i)}_j+\lambda \ell^{(i)}_j\}$
with a generic collection $\{\ell^{(i)}_j\}\in {\mathcal L}^{\bf k}$.
\end{proposition}

Let $\{\omega^{(i)}_j\}$ be a collection of germs of 1-forms on $(V, 0)$
with an isolated special point at the origin. Let $\nu : \widehat{V} \to V$
be the Nash transformation of the variety $V\subset B_\eps$ (see
Section~\ref{sect:Eu}). The collection of 1-forms $\{\omega^{(i)}_j\}$ gives
rise to a section $\widehat{\omega}$ of the bundle
$$
\widehat\TT=\bigoplus_{i=1}^s\bigoplus_{j=1}^{n-k_i+1}\widehat T^*_{i,j}
$$
where $\widehat{T}^\ast_{i,j}$ are copies of the dual Nash bundle
$\widehat{T}^\ast$ over the Nash transform $\widehat{V}$ numbered by
indices $i$ and $j$. Let $\widehat\DD\subset\widehat\TT$ be the set of
pairs $(x,\{\alpha^{(i)}_j\})$ where $x\in\widehat V$ and the collection
$\{\alpha^{(i)}_j\}$ of elements of $\widehat T_x^*$ (i.e.\ of linear
functions on $\widehat T_x$) is such that $\alpha^{(i)}_1$, \dots,
$\alpha^{(i)}_{n-k_i+1}$ are linearly dependent for each $i=1, \dots, s$.
The image of the section $\widehat\omega$ does not intersect $\widehat\DD$
outside of the preimage $\nu^{-1}(0)\subset\widehat V$ of the origin.
The map $\widehat\TT\setminus \widehat\DD\to \widehat V$ is a fibre
bundle. The fibre $W_x=\widehat\TT\setminus \widehat\DD$ of it is
$(2n-2)$-connected, its homology group $H_{2n-1}(W_x;\ZZ)$ is isomorphic
to $\ZZ$ and has a natural generator (see above). The latter fact implies
that the fibre bundle $\widehat\TT\setminus \widehat\DD\to \widehat V$
is homotopically simple in dimension $2n-1$, i.e.\ the fundamental group
$\pi_1(\widehat V)$ of the base acts trivially on the homotopy group
$\pi_{2n-1}(W_x)$ of the fibre, the last one being isomorphic to the
homology group $H_{2n-1}(W_x)$: see, e.g., \cite{Steenrod51}.
The following definition was made in \cite{EGBrasselet} (see also \cite{EGTrieste}).

\begin{definition}
The {\em local Chern obstruction} \index{local Chern obstruction} ${\rm Ch}(\{\omega^{(i)}_j\};V,0)$ of
the collections of germs of 1-forms $\{\omega^{(i)}_j\}$ on $(V,0)$ at the
origin is the (primary) obstruction to extend the section $\widehat{\omega}$
of the fibre bundle $\widehat\TT\setminus \widehat\DD\to \widehat V$ from
the preimage of a neighbourhood of the sphere $S_\eps= \partial B_\eps$ to
$\widehat V$, more precisely its value $($as an element of the homology group
$H^{2n}(\nu^{-1}(V\cap B_\eps), \nu^{-1}(V \cap S_\eps);\ZZ)$\,$)$ on the
fundamental class of the pair
$(\nu^{-1}(V\cap B_\eps), \nu^{-1}(V \cap S_\eps))$.
\end{definition}

The local Chern obstruction
${\rm Ch}(\{\omega^{(i)}_j\};V,0)$ can also be described as an intersection number, see \cite{EGBrasselet}. 
Namely, let ${\mathcal D}^{\bf k}_V\subset \CC^N \times {\mathcal L}^{\bf k}$
be the closure of the set of pairs $(x, \{\ell^{(i)}_j\})$ such that
$x \in V_{\rm reg}$ and the restrictions of the linear functions
$\ell_1^{(i)}$, \dots , $\ell_{n-k_i+1}^{(i)}$ to
$T_x V_{\rm reg} \subset \CC^N$ are linearly dependent for each
$i=1, \ldots, s$. (For $s=1$, ${\bf k}=\{n\}$, ${\mathcal D}^{\bf k}_V$ is
the (non-projectivized) conormal space of $V$ \cite{Teissier82}.) The
collection $\{\omega^{(i)}_j\}$ of germs of 1-forms on $(\CC^N,0)$
defines a section $\check{\omega}$ of the {trivial} fibre bundle
$\CC^N \times {\mathcal L}^{\bf k} \to \CC^N$. Then
\begin{equation} \label{eq:Chern_intersection}
{\rm Ch}(\{\omega^{(i)}_j\};V,0) = ( \check{\omega}(\CC^N) \circ
{\mathcal D}^{\bf k}_V)_0 ,
\end{equation}
where $(\cdot \circ \cdot )_0$ is the intersection number at the origin
in $\CC^N \times {\mathcal L}^{\bf k}$. This description can be considered as
a generalization of an expression of the local Euler obstruction as a
micro-local intersection number defined in \cite{KS90}, see also
\cite[Sections 5.0.3 and 5.2.1]{Schurmann03} and \cite{Schurmann04}.

\begin{remark}
 The local Euler obstruction is defined for vector fields on singular varieties as well as for 1-forms. One can see that collections of vector fields are not well adapted to a definition of the local Chern obstructions, at least on varieties with non-isolated singularities. The reason is as follows. Vector fields on a singular variety $V$ are required to be tangent to the smooth strata of $V$ (of a Whitney stratitfication). For example, on a one-dimensional stratum, all vector fields are proportional to each other and therefore a collection cannot have an isolated special point in a natural sense. A natural definition of the Chern obstruction of a collection of vector fields makes sense for varieties with isolated singularities. (Besides that, on a singular variety (with non-isolated singularities) continuous vector fields with isolated singular points exist, whereas this may not hold for holomorphic ones, see Remark~\ref{rmk:holomorphic_vf}.)
\end{remark}

Being a (primary) obstruction, the local Chern obstruction satisfies the
law of conservation of number, i.e.\ if a collection of 1-forms
$\{\widetilde\omega^{(i)}_j\}$ is a deformation of the collection
$\{\omega^{(i)}_j\}$ and has isolated special points on $V$, then
$$
{\rm Ch}(\{\omega^{(i)}_j\};V,0) = \sum {\rm Ch}(\{\omega^{(i)}_j\};V,p)
$$
where the sum on the right hand side is over all special points $p$ of
the collection $\{\widetilde\omega^{(i)}_j\}$ on $V$ in a neighbourhood
of the origin. With Proposition~\ref{prop:Morse} this implies the following
statements. The first statement is an analogue of  \cite[Proposition~2.3]{STV05b}.

\begin{proposition}\label{prop:no_special}
The local Chern obstruction ${\rm Ch}(\{\omega^{(i)}_j\};V,0)$ of a
collection $\{\omega^{(i)}_j\}$ of germs of holomorphic 1-forms is equal
to the number of special points on $V$ of a generic (holomorphic)
deformation of the collection (lying on $V_{\rm reg}$).
\end{proposition}

\begin{proposition}\label{prop:Chern=inv}
Let  $\{\omega^{(i)}_j\}$ be a collection of 1-forms on a compact (say,
projective) variety $V$ with only isolated special points. Then the sum
of the local Chern obstructions of the collection $\{\omega^{(i)}_j\}$ at
these points does not depend on the collection and therefore is an
invariant of the variety.
\end{proposition}

It is reasonable to consider this sum as ($(-1)^n$ times) the corresponding
Chern characteristic number of the singular variety $V$. It is well known that the characteristic numbers of a compact complex manifold cannot have arbitrary values, they satisfy certain divisibility properties. A.~Buryak \cite{BurMRL} showed that, contrary to this fact, any set of integers can be the set of Chern characteristic numbers of a singular projective variety.

Let $(V,0)$ be an {\icis}.
The fact that both the Chern obstruction and the GSV index of a collection
$\{\omega^{(i)}_j\}$ of 1-forms satisfy the law of conservation of number
and they coincide on a smooth manifold yields the following statement.

\begin{proposition}\label{prop4a}
For a collection $\{\omega^{(i)}_j\}$ on an {\icis} $(V,0)$, the difference
$$
\indGSV(\{\omega^{(i)}_j\};V,0)- {\rm Ch}(\{\omega^{(i)}_j\};V,0)
$$
does not depend on the collection and therefore is an invariant of the {\icis}.
\end{proposition}

In the framework of the definition of Schwartz Chern classes of singular varieties \cite{Schwartz65, BS81}, one has to consider $k$-fields on $n$-dimensional varieties. Let us give the definition of the Euler obstruction of a $k$-field following \cite{BSS09}.

Let $(V,0)\subset (\CC^N,0)$ be the germ of a pure $n$-dimensional complex analytic variety. Let $\CC^N= \bigcup_{i=1}^q V_i$ be a Whitney stratification of $\CC^N$ compatible with $V$, i.e.\ $\CC^N \setminus V$ is a stratum. Let (K) be triangulation of $\CC^N$ subordinated to the stratification  $\CC^N=\bigcup_{i=1}^q V_i$ and (D) be a cell decomposition of $\CC^N$ dual to (K).

Let $\{ X_j \}=(X_1 \ldots , X_k)$ be a $k$-field, i.e. a collection of stratified vector fields, i.e.\ at each point $p \in V$ each vector field $X_j$, $j=1, \ldots , k$, is tangent to the  stratum containing $p$. Assume that $\{ X_j \}$ has only isolated singular points in $V$. Let $\sigma$ be a $2(N-k+1)$-cell of (D). Note that $\sigma$ is transverse to all the strata $V_i$, $i=1, \ldots , q$. We assume that $\{ X_j \}$ has an isolated singularity at the barycentre $b_\sigma$ of $\sigma$ and is a $k$-frame in $(\sigma \setminus \{ b_\sigma \}) \cap V$, in particular it does not have any singularities on $\partial \sigma \cap V$. 

Let $\nu: \widehat{V} \to V$ be the Nash transform of $V$ and $\widehat{T}$ be the Nash bundle, (cf.\ Section~\ref{sect:Eu}). Each vector field $X_j$ lifts to a section $\widehat{X}_j$ of the bundle $\widehat{T}$ over $\nu^{-1}(\partial \sigma \cap V)$, see Section~\ref{sect:Eu}.  The $k$-frame $\{ X_j\}$ lifts to $k$ linearly independent sections  $\{ \widehat{X}_j \}$ of $\widehat{T}$ over $\nu^{-1}(\partial \sigma \cap V)$.

\begin{definition} The {\em local Euler obstruction}\index{local Euler obstruction! of a $k$-field}  ${\rm Eu}(\{ X_j \};V,\sigma)$  of the $k$-field $\{ X_j \}$ at $b_\sigma$ is the obstruction to extend $\{ \widehat{X}_j \}$ to a collection of $k$ linearly independent sections of $\widehat{T}$ over $\nu^{-1}(\sigma \cap V)$, more precisely its value (as an element of $H^{2(N-k+1)}(\nu^{-1}(\sigma \cap V),\nu^{-1}(\partial \sigma \cap V))$) on the fundamental class of the pair $(\nu^{-1}(\sigma \cap V),\nu^{-1}(\partial \sigma \cap V))$.
\end{definition}

Let $\{ \omega_j \}$ be a collection of 1-forms on $V$ with an isolated singularity at the barycentre $b_\sigma$ of $\sigma$.

\begin{definition} The {\em local Euler obstruction} ${\rm Eu}(\{ \omega_j \}; V, \sigma)$ of the collection $\{ \omega_j \}$ at $b_\sigma$ is defined in a similar way, but now taking sections of the dual Nash bundle $\widehat{T}^\ast$.
\end{definition}

There is also the definition of an Euler obstruction of a map due to N.~Grulha \cite{Grulha08}.

Let $f=(f_1, \ldots , f_k): (V,0) \to (\CC^k,0)$ be the germ of an analytic map. Let ${\rm grad}_V f_j$, $j=1, \ldots , k$, be the vector fields constructed in Section~\ref{sect:Eu}. The construction can be done in such a way that for $x \in V \setminus \{ 0 \}$ the vector fields $({\rm grad}_V f_1(x), \ldots , {\rm grad}_V f_k(x))$ are linearly independent, see \cite{Grulha08}.

Let $\Sigma f$ be the singular set of $f$. Let us assume that there exists a cell decomposition (D) of $\CC^N$ and a $2(N-k+1)$-cell of (D) with barycentre 0 such that  
$\Sigma f \cap \partial \sigma=\emptyset$.

\begin{definition}
The {\em local Euler obstruction }\index{local Euler obstruction! of a map} ${\rm Eu}(f;V,\sigma)$ of the map $f$ relative to $\sigma$ is defined to be
\[ {\rm Eu}(f;V,\sigma):= {\rm Eu}(\{ {\rm grad}_V f_j \};V,\sigma).
\]
\end{definition}

It is shown in \cite[Corollary~5]{BGR10} that the definition of ${\rm Eu}(f;V,\sigma)$ does not depend on a generic choice of the cell $\sigma$. 

Instead of the collection of vector fields $\{ {\rm grad}_V f_j\}$, we can consider the collection $\{ df_j \}$ of 1-forms associated to $f$. This leads to the following definition (see \cite{BGR10}).

\begin{definition}
\[ {\rm Eu}^\ast(f;V,\sigma):= {\rm Eu}(\{ df_j \};V, \sigma).
\]
\end{definition}

Relations between the Euler obstructions of $k$-fields and Chern obstructions were described in \cite{BGR10}.

Let $(V,0)\subset (\CC^N,0)$, $\{ V_i \}$, (K), and (D) be as above. Assume that 0 is the barycentre of a $2k$-simplex of the triangulation (K) and let $\sigma$ be the dual $2(N-k)$-cell. Since (K) is subordinated to the stratiication, the simplex $\tau$ is contained in a stratum and the cell $\sigma$ is tranverse to the strata. A neighbourhood of 0 in $\CC^N$ is homeomorphic to $\sigma \times \tau$ and one has 
\[ (\sigma \times \tau) \cap V= (\sigma \cap V) \times \tau \cong V \cap B_\eps
\]
for a ball $B_\eps$ around the origin of sufficiently small radius $\eps$. Denote by $\pi_1$ and $\pi_2$ the projections $\pi_1: \sigma \times \tau \to \sigma$ and $\pi_2: \sigma \times \tau \to \tau$. Note that $\tau$ is a smooth manifold, so the index  $\ind(\{ \omega_j \};\tau,0)$ of a collection of 1-forms $\{ \omega_j \}$ according to Section~\ref{sect:coll} is defined. The following theorem is proved in \cite[Theorem~2.2]{BGR10}.

\begin{theorem}[Brasselet, Grulha, Ruas]
In the above setting, let $\{ \omega_j^{(1)} \}$, $j=1, \ldots , k-1$, be a collection of germs of 1-forms on $\sigma$ and $\{ \omega_j^{(2)} \}$, $j=1, \ldots , d-k+1$, be a collection of germs of 1-forms on $\tau$. The collection of germs of 1-forms on $(\CC^N,0)$ given by $\{ \omega_j^{(i)} \} = \{ \pi_1^\ast(\omega_{j_1}^{(1)}),\pi_2^\ast(\omega_{j_2}^{(1)}) \}$ satisfies
\[
{\rm Ch}(\{\omega^{(i)}_j\};V,0) = {\rm Eu}(\{ \omega_j^{(1)} \};V,\sigma) \cdot \ind(\{ \omega_j^{(2)} \};\tau,0).
\]
\end{theorem}

The following corollary is derived from this theorem, see \cite[Corollary~2.3, Corollary~2.6]{BGR10}.

\begin{corollary}[Brasselet, Grulha, Ruas]
Let $(V,0)$ be as above and let $f :(V,0) \to (\CC^k,0)$ be a map germ. Let $\{ \omega_j^{(i)} \}=\{ \omega_{j_1}^{(1)}, \omega_{j_2}^{(i)} \}$ be the collection of 1-forms defined by $\{ \omega_{j_1}^{(1)} \} = \{ df_1, \ldots , df_k \}$ and $\{ \omega_{j_2}^{(2)} \} = \{ \ell_1, \ldots , \ell_{n-k+2} \}$ where $\ell_1, \ldots, \ell_{n-k+1}$ are linearly independent linear forms dual to the tangent field of $\sigma$ and $\ell_{n-k+2}$  is a radial linear form. Then one has
\begin{equation*}
{\rm Eu}^\ast(f;V, \sigma) = {\rm Ch}(\{\omega^{(i)}_j\};V,0) = ( \check{\omega}(\CC^N) \circ
{\mathcal D}^{\bf k}_V)_0 \, 
\end{equation*}
(cf.\ Equation~\ref{eq:Chern_intersection}).
\end{corollary}

It follows from this corollary that ${\rm Eu}^\ast(f;V, \sigma)$ is independent of a generic choice of $\sigma$. Moreover, the following identity \cite[Theorem~2.4]{BGR10} is proved using this corollary:
\[ {\rm Eu}(f;V, \sigma) = (-1)^{n-k+1} {\rm Eu}^\ast(f;V, \sigma).
\]
Therefore the Euler obstruction ${\rm Eu}(f;V, \sigma)$ is also independent of a generic choice of $\sigma$ \cite[Corollary~2.5]{BGR10}.

In \cite{BGR10}, also some formulas to compute the Chern obstruction are given. In particular, it is shown that the Chern obstruction is related with the polar multiplicity. \index{polar multiplicity}\index{multiplicity! polar}
Let $(V,0) \subset (\CC^N,0)$ be the germ of a pure $n$-dimensional complex analytic variety and $f:(V,0) \to \CC^k$ be a generic projection. The $(n-k+1)$-polar variety $P_{n-k+1}(V)$ is the closure of the singular set $\overline{\Sigma f}$ of $f$. Its multiplicity at 0 is denoted by $m_{n-k+1}(V,0)$. The following theorem is proved  in \cite[Theorem~3.1]{BGR10}.

\begin{theorem}[Brasselet, Grulha, Ruas]
Let $(V,0) \subset (\CC^N,0)$ be the germ of a pure $n$-dimensional complex analytic variety and $f:(V,0) \to \CC^k$ be a generic projection. Let $\{\omega^{(i)}_j\}$ be a collection of germs of 1-forms on $(\CC^N,0)$ such that $\{\omega^{(1)}_j\}=\{ df_1, \ldots , df_k \}$ and $\{\omega^{(i)}_j\}$, $i=2, \ldots ,s$, $j=1, \ldots , n-k_i+1$, are generic subcollections, where $k_i$, $i=2, \ldots, s$, are non-negative integers with $\sum_{i=2}^s k_i=k-1$. Then one has
\[ {\rm Ch}(\{\omega^{(i)}_j\};V,0) = m_{n-k+1}(V,0).
\]
\end{theorem}

In \cite[Theorem~6.1]{GGr}, a formula for the Chern obstruction of a collection of 1-forms on an equidimensional analytic variety is given in terms of the multiplicity of pairs of modules.

\subsection{Homological index} \label{sect:Homcoll}
Here we consider the definition of the homological index for a collection of 1-forms due to E.~Gorsky and the second author \cite{GGZ}.

Let $(V,0)$ be the germ of a complex algebraic variety of dimension $n$ with an isolated singular point at the origin. Let $k_i$, $i=1, \ldots ,s$, be positive integers such that $\sum_{s=1}^s k_i=n$ and let $\{\omega^{(i)}_j\}$, $i=1, \ldots, s$, $j=1, \ldots , n-k_i+1$, be a collection of germs of holomorphic 1-forms on $(V,0)$. 

Let $W_i=\CC^{n-k_i+1}$ be an auxiliary vector space with a basis $u_1, \ldots , u_{n-k_i+1}$. We consider the complex $\calC^{(i)} = \calC(\omega^{(i)}_1, \ldots , \omega^{(i)}_{n-k_i+1})$ of sheaves of $\calO_{V,0}$-modules defined as follows:
\[
\calC_0^{(i)}:=\Omega^n_{V,0}, \quad \calC_t^{(i)}:=\Omega^{k_i-t} \otimes S^{t-1}W_i, \quad 1 \leq t \leq k_i.
\]
The differential $d_t:\calC_t^{(i)} \to \calC_{t-1}^{(i)}$ is defined by
\begin{eqnarray*}
d_1(\beta) & := & \beta \wedge \omega_1^{(i)} \wedge \ldots \wedge \omega_{n-k_i+1}^{(i)}, \\
d_t(\beta \otimes \varphi(u)) &:=& \sum_{l=1}^{n-k_i+1} \left( \beta \wedge \omega_l^{(i)} \right) \otimes \frac{\partial \varphi}{\partial u_l}, \quad 2 \leq t \leq k_i.
\end{eqnarray*}
The complex $(\calC^{(i)}, d)$ is indeed a chain complex \cite[Lemma~11]{GGZ}. It is shown in \cite[Lemma~14]{GGZ} that the cohomology groups of $\calC^{(i)}$ are supported on the locus of the points where the forms $\{\omega^{(i)}_j\}$ are linearly dependent. Define
\[ \calC= \bigotimes_{i=1}^s \calC^{(i)},
\]
where the tensor product is taken over $\calO_{V,0}$.

\begin{definition}
The {\em homological index}\index{homological index! of a collection} $\indhom(\{\omega^{(i)}_j\};V,0)$  of the collection of 1-forms $\{\omega^{(i)}_j\}$ is the Euler characteristic of the complex $\calC$:
\[
\indhom(\{\omega^{(i)}_j\};V,0):= \sum_{t=0}^n (-1)^t \dim H^t(\calC).
\]
\end{definition}

The homological index for a collection of 1-forms with an isolated singular point satisfies the law of conservation of number \cite[Proposition~18]{GGZ}. The following theorem is \cite[Theorem~19]{GGZ}.

\begin{theorem}[Gorsky, Gusein-Zade]
Let $(V,0)$ be an {\icis} and let $\{\omega^{(i)}_j\}$ be a collection of holomorphic 1-forms on $(V,0)$ with an isolated singular point. Then
\[ \indhom(\{\omega^{(i)}_j\};V,0)=\indGSV(\{\omega^{(i)}_j\};V,0).
\]
\end{theorem}

Since both the homological index and the Chern obstruction satisfy the law of conservation of number and coincide on a smooth manifold, one has the following statement (see \cite[Proposition~40]{GGZ}).

\begin{proposition}
Let $(V,0) \subset (\CC^N,0)$ be the germ of a complex analytic variety of pure dimension $n$ with an isolated singularity at the origin. The difference
\[
\indhom(\{\omega^{(i)}_j\};V,0) - {\rm Ch}(\{\omega^{(i)}_j\};V,0)
\]
between the homological index and the Chern obstruction does not depend on the collection $\{\omega^{(i)}_j\}$ and is an invariant of the singularity $(V,0)$.
\end{proposition}


\section{Equivariant indices}\label{sec:equiv_ind}

\subsection{Equivariant Euler characteristics}\label{subsec:equiv_Euler}
The notions of indices of vector fields and of 1-forms (on smooth manifolds and on singular varieties) are related with the Euler characteristic (through the Poincar\'e--Hopf theorem). Therefore it is natural to discuss equivariant versions of the Euler characteristic first.

In what follows we shall consider the additive Euler characteristic defined (for topological spaces nice enough, say, for those homeomorphic to
locally compact unions of cells in finite CW-complexes) as the alternating sum of the
dimensions of the cohomology groups with compact support:
\begin{equation}\label{chi_add}
 \chi(V)=\sum_{q=0}^{\infty}(-1)^q \dim H^q_{\rm c}(V;\CC)\,.
\end{equation}
This Euler characteristic coincides with the ``traditional one'' (defined as the alternating sum of the dimensions of the usual cohomology groups) for compact spaces (finite CW-complexes) and for complex quasi-projective varieties. The additivity of the
Euler characteristic permits to use it as
a sort of a (non-positive) measure for the
definition of the integral with respect to
the Euler characteristic: \cite{Viro}.

There are several generalizations of the notion of the Euler characteristic to the equivariant setting, i.e.\ for
spaces with actions of a group (say, a finite one).
The most simple (and the most straightforward) one is obtained by substituting the dimensions of the cohomology groups in~(\ref{chi_add}) by the classes of the corresponding
$G$-modules $H^q_{\rm c}(V;\CC)$
(spaces of representations of $G$) in the ring $R(G)$ of representations of the group $G$.
This analogue of the Euler characteristic is defined as an element of the ring $R(G)$.
It was introduced in \cite{Verdier} and was used, e.g., in \cite{CTCWall}.

A finer
equivariant version of the Euler characteristic of a $G$-space can be defined as an element of the Burnside ring $A(G)$ of the group $G$, i.e.\ the Grothendieck ring of finite $G$-sets.
The latter is the abelian group generated by the classes $[(Z,G)]$ of finite $G$-sets
modulo the following relations:
\begin{itemize}
 \item[---] if $(Z_1,G)$ and $(Z_2,G)$ are isomorphic, i.e., if there exists a bijective
 $G$-equivariant map $Z_1\to Z_2$, then $[(Z_1,G)]=[(Z_2,G)]$;
 \item[---] $[(Z_1\sqcup Z_2,G)]=[(Z_1,G)]+[(Z_2,G)]$.
\end{itemize}
The multiplication in $A(G)$ is defined by the Cartesian product of sets with the 
natural (diagonal) $G$-action. The Burnside ring $A(G)$ is the free abelian group generated by the classes of
irreducible $G$-sets which are in bijection with the classes of conjugate
subgroups of the group $G$: the conjugacy class $[H]$ of a subgroup $H\subset G$ corresponds to the class
of the $G$-set $G/H$.
One has a natural ring homomorphism $A(G)\to R(G)$, sending a $G$-set $Z$ to the space
of functions on $Z$ with the natural (left) action of the group $G$: $(g^*f)(y)=f(g^{-1}y)$. 
This homomorphism is, in general, neither a monomorphism nor an epimorphism.
In what follows we shall mostly use the
equivariant version of the Euler characteristic with values in $A(G)$ and therefore we shall refer to it as the {\em equivariant Euler characteristic}.

Let $V$ be a sufficiently nice space with an action of the group $G$. 
For a point $x$ of the space $V$ let $G_x$ be the isotropy
subgroup $\{g\in G:gx=x\}$ of the point $x$. For a subgroup $H$ of the group $G$ let $V^H$
be the fixed point set of the group $H$: 
$\{x\in V: G_x\supset H\}$, and let $V^{(H)}$
be the set $\{x\in V: G_x= H\}$ of points with the isotropy subgroup coinciding with $H$.
Let ${\rm Conjsub\,}G$ be the set of the conjugacy classes of subgroups of $G$.
For a conjugacy class $[H]$ of subgroups of $G$
(which contains the subgroup $H$), let $V^{[H]}$ be the set of points such that each of them
is fixed with respect to a subgroup conjugate to $H$ and let $V^{([H])}$ be the set of points $x\in V$ whose isotropy subgroups $G_x$ are conjugate to $H$. 

\begin{definition}
 The {\em equivariant Euler characteristic}\index{equivariant Euler characteristic}\index{Euler characteristic! equivariant}
 of the $G$-space $V$ is defined by 
 \begin{equation}\label{EquiEulerGen}
 \chi^G(V)=\sum_{[H]\in{\rm Conjsub\,}G}\chi(V^{([H])}/G)\cdot[G/H]\in A(G)\,.
 \end{equation}
\end{definition}

This notion was introduced in \cite{TtD}.
The equivariant Euler characteristic satisfies the additivity property:
if $W$ is a closed $G$-invariant subspace of a $G$-space $V$, then
$$
\chi^G(V)=\chi^G(W)+\chi^G(V\setminus W)\,.
$$
One can show that it is a universal invariant
possessing this property (on the class of spaces homeomorphic to locally closed unions of cells in finite CW-complexes with cell
actions of the group $G$).
The equivariant analogue of the Euler characteristic with values in the ring $R(G)$
of representations of $G$ is obtained from
the equivariant Euler characteristic by
the natural homomorphism $A(G)\to R(G)$
described above.
Among other reductions of the equivariant Euler characteristic, one can indicate the
orbifold Euler characteristic (see, e.g.,
\cite{AS}, \cite{HH}) and its higher order analogues (\cite{AS}, \cite{Tamanoi}).

\subsection{Equivariant indices of vector fields and 1-forms on manifolds}\label{subsec:equiv_on_smooth}
Let a finite group $G$ act smoothly on (the germ of) the affine
space $(\RR^n,0)$. Without loss of generality one can assume that
the action is linear, i.e.\ it is defined by a representation of $G$ on $\RR^n$. Let $X$ be a (continuous) vector field on $(\RR^n,0)$ invariant with respect to the action of $G$ and with an isolated singular point at the origin. (One can see that, for each
point $p$ from a neighbourhood of the origin (where $X$ is defined), the vector $X(p)$ is tangent to the subspace
$(\RR^n)^{G_p}$ of the fixed points of the isotropy subgroup
$G_p$ of the point $p$.)

We assume $\RR^n$ to be endowed with a $G$-invariant Euclidean metric. Let $\eps>0$ be small enough
so that the vector field $X$ is defined on a neighbourhood of
the closed ball $B_\eps$ of radius $\eps$ centred at the origin and 
has no singular points in $B_\eps$ outside the origin. It is easy to see that there
exists a $G$-invariant vector field $\widetilde{X}$ on a neighbourhood of $B_\eps$ such that:
\begin{itemize}
\item[(1)] The vector field $\widetilde{X}$ coincides with $X$ on a neighbourhood of the sphere $S_\eps= \partial B_\eps$.
\item[(2)] In a
neighbourhood of each singular point $p\in B_\eps\setminus\{0\}$, the vector field $\widetilde{X}$ is as follows.
Let $H=G_{p}$ be the isotropy subgroup of the point $p$. The germ $(\RR^n,p)$ is in a natural way isomorphic to 
$((\RR^n)^H,p)\times (((\RR^n)^H)^{\perp},0)$,
where $((\RR^n)^H)^{\perp}$ is the orthogonal complement to the
subspace $(\RR^n)^H$: the direct sum of the subspaces of $\RR^n$
corresponding to non-trivial representations of $H$.
In a neighbourhood of $p$ the vector $\widetilde{X}(y_1,y_2)$
($y_1\in ((\RR^n)^H,x_0)$, $y_2\in ((\RR^n)^H)^{\perp},0)$)
is the sum $X_1(y_1)+X_2(y_2)$, where
$X_1$ is a vector field on $((\RR^n)^H,p)$ with an isolated
singular point at $p$, $X_2$ is an $H$-invariant radial vector
field on $((\RR^n)^H)^{\perp},0)$.
(Let us recall that, on a zero-dimensional space
the only vector field (zero) is non-degenerate with the index 1 and also radial.)
\end{itemize}

\begin{remark} One can assume that the vector field $X_1$ is smooth and has a non-degenerate
singular point at $x_0$ (and therefore $\ind(X_1; (\RR^n)^H,x_0) = \pm 1$), however, this is not necessary for the definition.
\end{remark} 

\begin{definition} The {\em equivariant index}\index{equivariant index}\index{index! equivariant} 
 $\ind^G(X; \RR^n,0)$ of the vector field $X$ at the origin
 is defined by the equation 
 \begin{equation*}\label{Defrad}
\ind^G(X; \RR^n,0)=\sum_{\overline{p}\in ({\Sing}\widetilde{X})/G} 
\ind(\widetilde{X}_{\vert V_{(p)}}; V_{(p)},p) [Gp]\,, 
 \end{equation*}
 where $p$ is a representative of the orbit $\overline{p}$. (One has $[Gp]=[G/G_p]$.)
 \end{definition}

\begin{remark}
One can say that the equivariant index $\ind^G(X; \RR^n,0)$ is the class $[\Sing \widetilde{X}]\in B(G)$ of the set $\Sing \widetilde{X}$ of singular
 points of $\widetilde{X}$ with multiplicities equal to
the usual indices $\ind(\widetilde{X}_{\vert (\RR^n)^{G_p}}; (\RR^n)^{G_p},p)$ of the restrictions
of the vector field $\widetilde{X}$ to the corresponding 
fixed point sets (smooth manifolds).
\end{remark}

One has the following equivariant version of the Poincar\'e-Hopf theorem.

For a subgroup $H\subset G$ there are natural maps
$\mbox{R}_{H}^{G}: B(G)\to B(H)$ and $\mbox{I}_{H}^{G}: B(H)\to B(G)$.
The {\em restriction map} $\mbox{R}_{H}^{G}$ sends a $G$-set $Z$ to the same set considered with the $H$-action.
The {\em induction map} $\mbox{I}_{H}^{G}$ sends an $H$-set $Z$ to the product $G\times Z$ factorized
by the natural equivalence: $(g_1,
 x_1)\sim (g_2, x_2)$ if there exists $g\in H$ such that
$g_2=g_1g$, $x_2=g^{-1}x_1$ with the natural (left) $G$-action. Both maps are group homomorphisms,
however the induction map $\mbox{I}_{H}^{G}$ is not a ring homomorphism. 

\begin{theorem} \label{thm:PHMequ}
Let $M$ be a closed (compact, without boundary)
$G$-manifold and let $X$ be a $G$-invariant vector field on $M$ with isolated singular
points. Then one has
\begin{equation*}
 \sum_{\overline{p}\in (\Sing X)/G} 
{\rm I}_{G_p}^G (\ind^{G_p}(X; M,p))=\chi^G(M)\,.
\end{equation*}
\end{theorem}

A version of this definition of an equivariant index was given first in~\cite{Luck}. However,
the index there takes values in an extension
of the Burnside ring which depends on the $G$-manifold $M$ and is not a ring. (It takes into account connected components of the fixed point sets of subgroups of $G$.)

Almost the same definition can be given for the equivariant index $\ind^{G}(\omega; \RR^n,0)$
of a $G$-invariant 1-form $\omega$ on $(\RR^n,0)$. A $G$-invariant Riemannian metric on a $G$-manifold permits to identify invariant 1-forms with invariant vector fields. In this way,
the equivariant index of a 1-form is the equivariant index of the corresponding vector field.

In the complex setting one has the usual sign correction factor for a complex-valued 1-form on a complex manifold (see Section~\ref{subsec:smooth_complex}).

\subsection{The equivariant radial index on a singular variety}\label{subsec:equiv_radial_on_singular}
On a singular $G$-variety $G$-invariant vector
fields and $G$-invariant 1-forms cannot be identified with each other. Therefore their equivariant indices (even being defined in similar ways) cannot be expressed through each other.

Here we shall give the definition of an
equivariant version of the radial index
for vector fields. The necessary changes for 1-forms are clear.

Let $(V,0)$ be the germ of a closed real subanalytic set with an action
of a finite group $G$. We assume $(V,0)$ to be embedded into $(\RR^N, 0)$
and the $G$-action to be induced by an (analytic) action on a neighbourhood of
the origin in $(\RR^N, 0)$. (The action on $(\RR^N, 0)$ can be assumed to be linear.) 

Let $V=\bigcup_{i=1}^q V_i$ be a subanalytic Whitney $G$-stratification of $V$. This means that each stratum $V_i$ is $G$-invariant,
the isotropy subgroups $G_p=\{g\in G: gp=p\}$ of all points $p$ of $V_i$ are conjugate to each other,
and the quotient of the stratum $V_i$ by the group $G$
is connected.

Let $X$ be a $G$-invariant (stratified) vector field on $(V,0)$ with an isolated singular point
at the origin. One can show that there exists a (continuous) $G$-invariant
stratified vector field $\widetilde{X}$ on $V$ satisfying Conditions (1)--(3) from Section~\ref{sect:rad}.

The following definition was made in \cite{EGEJM}.
Let $A$ be the set (a $G$-set) of the singular points of the vector field
$\widetilde{X}$ on $V\cap B_\eps$ considered with the multiplicities equal to
the usual indices $\ind(\widetilde{X}_{\vert V_{(p)}}; V_{(p)},p)$ of the restrictions
of the vector field $\widetilde{X}$ to the corresponding strata (smooth manifolds).

\begin{definition}
 The {\em equivariant radial index}\index{equivariant radial index}\index{radial index! equivariant} $\indrad^G(X; V,0)$ of the vector field $X$
 on $V$ at the origin is the class $[A]\in B(G)$ of the set $A$ of singular
 points of $\widetilde{X}$ with multiplicities.
\end{definition}

One can show that the equivariant radial index is well-defined: see~\cite{EGEJM}.

\begin{remark}
 As above (in the smooth case) one can write the definition as
\begin{equation*}
\indrad^G(X; V,0)=\sum_{\overline{p}\in ({\Sing}\widetilde{X})/G} 
\ind(\widetilde{X}_{\vert V_{(p)}}; V_{(p)},p) [Gp]\,, 
\end{equation*}
 where $p$ is a representative of the orbit $\overline{p}$.
\end{remark}

For a subgroup $H\subset G$, the vector field $X$ is $H$-invariant and one has
$\indrad^H(X; V,0)=
\mbox{R}^G_H(\indrad^G(X; V,0))$.

One has the following generalization of Theorem~\ref{thm:PHMequ} (see \cite[Theorem~4.6]{EGEJM}).
\begin{theorem} 
Let $V=\bigcup_{i=1}^q V_i$ be a compact subanalytic variety and  let $X$ be a $G$-invariant stratified vector field on $V$ with isolated singular
points. Then one has
\begin{equation*}
 \sum_{\overline{p}\in (\Sing X)/G} 
{\rm I}_{G_p}^G (\ind^{G_p}(X; V,p))=\chi^G(V)\,.
\end{equation*}
\end{theorem}

Another analogue of the Euler characteristic can be defined for orbifolds: the universal Euler characteristic of orbifolds~\cite{GLM18}. It takes values in the ring $\calR$ generated, as a free abelian group, by isomorphism classes of finite groups. The corresponding analogue of the radial
index of vector fields and 1-forms (with values in the same ring $\calR$) was defined in~\cite{Gusein-Zade21}.

\subsection{Equivariant GSV and Poincar\' e--Hopf index}\label{subsec:equiv_GSV}
Let the space $\CC^{N}$ be endowed with an action of a finite group $G$ (say, with a linear one) and 
let $(V,0)=\{ z \in(\CC^{N},0): f_1(z)=\cdots =f_{N-n}(z)=0\}$ be an $n$-dimensional germ of 
an isolated 
complete intersection singularity defined by $G$-invariant function germs $f_i:(\CC^{N},0)\to(\CC,0)$, $i=1,\ldots, {N-n}$.
In the usual (non-equivariant) setting (thus for the trivial group $G$), the GSV index of a vector field or of a 1-form on $(V,0)$ can be defined in terms of the degree of a certain map or in terms of the intersection number of some cycles. Equivariant versions of these notions (the degree and the intersection index) are not defined (at least as elements of the Burnside ring $A(G)$).
Therefore, in order to define equivariant versions of them, one has to use 
the fact that the GSV index agrees with the Poincar\'e--Hopf index (Proposition~\ref{prop:PH=GSV}). Namely, let $X$ be a $G$-invariant 
vector field on $(V,0)$ with an isolated singularity at the origin. Let $F:(\calV,0) \to (\CC,0)$ be the essentially unique smoothing of $(V,0)$ (see Section~\ref{sect:PH}). The vector field $\widetilde{X}$ of Section~\ref{sect:PH} can be made $G$-invariant. (To get a $G$-invariant vector field, one can take an arbitrary one and take the mean over the group.) Then one can make the following definition (cf.\ \cite[Definition~5.1]{EGEJM}).

\begin{definition} The {\em $G$-equivariant GSV index} (or {\em $G$-equivariant Poincar\' e--Hopf index}) is 
\begin{equation} \label{eq:GGSV}
\indGSV^G(X;V,0) = \sum_{\overline{p} \in ({\rm Sing}\, \widetilde{X})/G} \ind^G(\widetilde{X}; V_t^F, p) \in A(G),
\end{equation}
where $V_t^F$ is the Milnor fibre corresponding to the smoothing $F$. 
\end{definition}

It is easy to show that the right hand side of (\ref{eq:GGSV}) does not depend on the extension $\widetilde{X}$ and therefore the equivariant GSV index is well-defined.

There is also a generalization to non-isolated complete intersection singularities, see \cite[Definition~5.1]{EGEJM}.

The relation between the equivariant GSV index and the radial one can be described as follows. The Milnor fibre $V_t^F$ is a manifold with a $G$-action and therefore its equivariant Euler characteristic $\chi^G(V_t^F) \in A(G)$ is defined. Let $\overline{\chi}^G(V_t^F) = \chi^G(V_t^F) -1$ be the reduced equivariant Euler characteristic of $V_t^F$. (The element $(-1)^n\overline{\chi}^G(V_t^F) \in A(G)$ can be regarded as an equivariant version of the Milnor number of the {\icis} $(V,0)$). There is the following generalization of Proposition~\ref{prop:muX} (see \cite[Proposition~5.3]{EGEJM}).
\begin{proposition} 
$\indGSV^G(X; V,0) =\indrad^G(X;V,0)+\overline{\chi}^G(V_t^F)$.
\end{proposition}
 
 Let the group $G$ act on the projective space $\CC\PP^N$ by projective transformations and let $V \subset \CC\PP^N$ be a $G$-invariant complete intersection with isolated singularities. It has a natural $G$-invariant smoothing
 $\widetilde{V} \subset \CC\PP^N$. Let $X$ be a $G$-invariant vector field on $V$ with isolated singular points. Since the GSV index was defined as the Poincar\' e--Hopf index, it counts singular points of the vector field on the smoothing of the \icis. Therefore one has the following version of the Poincar\' e--Hopf theorem (see \cite[Proposition~5.2]{EGEJM}).
\begin{proposition}
Let $V$, $\widetilde{V}$,  and  $X$ be as above.  Then one has
$$
\sum_{\overline{p} \in(\Sing X)/G} {\rm I}_{G_p}^G( \indGSV^G(X;V,p)) = \chi^G(\widetilde{V})\in A(G)\,.
$$
\end{proposition}

\subsection{Equivariant homological index}\label{subsec:equiv_homological}
Let $\CC^N$ be endowed with an action (say, a linear one) of a finite group $G$ and let
$(V, 0)\subset (\CC^N, 0)$ be a $G$-invariant germ of an analytic variety of pure dimension $n$.
Let $X$ be a $G$-invariant holomorphic vector field on 
$(V,0)$ with an isolated singular point at the origin.

\begin{remark}
 The condition that such a vector field on $(V,0)$ exists is a rather restrictive condition on the variety. Namely, a neighbourhood in $V$ of any point $p$ of $V\setminus\{0\}$ has to be isomorphic to the direct product $(W_p,0)\times(\CC,0)$ for a variety $W_p$ (cf.\ Remark~\ref{rmk:holomorphic_vf}). This holds, in particular, if $V \setminus \{ 0 \}$ is non-singular.
\end{remark}

Let $\Omega^i_{V,0}$, $i=1,2,\ldots$, be the modules of germs of differential forms on $(V, 0)$
($\Omega^i_{V,0}=\calO_{V,0}$). One has natural actions of the group $G$ on them. Consider the complex~(\ref{eqn:hom_vector}).
It consists of $G$-modules and therefore its homology groups are (finite dimensional) $G$-modules as well.


\begin{definition}
 The {\em equivariant homological index}\index{equivariant homological index}\index{homological index! equivariant} of the vector field $X$ on $(V,0)$ is
 $$
 \indhom^G(X; V,0)=\sum_{i=0}^n
 \left[H_i(\Omega^{\bullet}_{V,0},X)\right]\in R_{\CC}(G)\,,
 $$
 where $R_{\CC}(G)$ is the ring of (complex) representations of the group $G$ and $[\cdot]$
 is the class of a $G$-module in $R_{\CC}(G)$.
\end{definition}

Let $\omega$ be a $G$-invariant holomorphic 1-form on 
$(V,0)$ with an isolated singular point at the origin.

\begin{remark}
 One can see that (in contrast to the situation for vector fields) 1-forms with this property always exist.
\end{remark}

Consider the complex~(\ref{eqn:hom_1-form}).
If $(V,0)$ is smooth and $\omega(0)\ne 0$, the homology groups of the complex (\ref{eqn:hom_1-form}) are trivial. 
This implies that, if both $V$ and $\omega$ have isolated singular points at the origin, the homology groups $H_i(\Omega^{\bullet}_{V,0},\wedge\omega)$ of the complex (\ref{eqn:hom_1-form}) are finite dimensional $G$-modules.

\begin{definition}
 Let $(V,0)$ have an isolated singular point at the origin and let $\omega$ have an isolated singular point at $0\in V$.
 The {\em equivariant homological index}\index{equivariant homological index}\index{homological index! equivariant} of the 1-form $\omega$ on $(V,0)$ is
 $$
 \indhom^G(\omega; V,0)=\sum_{i=0}^n
 \left[H_i(\Omega^{\bullet}_{V,0},\wedge\omega)\right]\in R_{\CC}(G)\,.
 $$
\end{definition}

\begin{remark}
 It is not clear whether this definition makes sense for an arbitrary ($G$-invariant) variety $(V, 0)\subset (\CC^N, 0)$, not necessarily with an isolated singular point at the origin, i.e.\ whether the homology groups
 $H_i(\Omega^{\bullet}_{V,0},\wedge\omega)$ are finite dimensional in this case as well.
\end{remark}

Assume that $(V,0)=(\CC^n,0)$. It is not difficult to show that the equivariant homological index of a vector field $X$ with an isolated singular point coincides with the reduction (under the natural homomorphism
$A(G)\to R_{\CC}(G)$) of the equivariant (radial) index of $X$. This follows from the fact that a $G$-invariant vector field on $\CC^n$ can be deformed to one with only non-degenerate singular points and it is easy to verify the coincidence for a vector field with a non-degenerate singular point. The situation is quite different for 1-forms (for a non-trivial group $G$).
A $G$-invariant 1-form cannot, in general, be deformed to one with non-degenerate singular points. It is not clear how one can describe non-removable singularities of invariant 1-forms for an arbitrary finite group $G$
(in order to verify the coincidence of the described indices for them). This was made in~\cite{Mamedova} for the group $\ZZ_3$ of order 3. For an arbitrary finite group $G$ the statement about the coincidence was proved in~\cite{GM17}. 

\subsection{Equivariant Euler obstruction}\label{subsec:equiv_Euler_obstruction}
The Euler obstruction of a vector field or of a 1-form
at an isolated singular point on a (quasi-projective) singular variety can be regarded as a version of an index
of it. Similar to the GSV index, the usual (non-equivariant) version of the Euler obstruction is defined in terms of the (first) obstruction to extend a section of a bundle. An equivariant version of the notion of the first obstruction is not defined (at least as an element of the Burnside ring $A(G)$). Therefore, to define an equivariant
version of the Euler obstruction (of a vector field or of a 1-form), one has to use another approach.

A method to define the (local) equivariant Euler obstruction of an invariant 1-form was suggested in~\cite{EGBrasil}. The idea resembles the one used for the definition of the equivariant radial index.
Let $(V,0)\subset(\CC^N,0)$ be a germ of a complex analytic variety with an action of a finite group $G$
and let $\{V_i\}_{i \in I}$  be a $G$-invariant Whitney stratification of it.
Let $\omega$ be a germ of a $G$-invariant complex
1-form on $(V,0)$ (that is the restriction of a $G$-invariant 1-form on $(\CC^N,0)$) with an isolated singular point at the origin. Let $B_\eps$ be a ball of a small radius
$\eps$ around the origin such that representatives of $V$ and of $\omega$ are defined in $B_\eps$
and the 1-form $\omega$ has no singular points on $V\setminus\{0\}$ inside $B_\eps$.
Let $\widetilde\omega$ be a $G$-invariant 1-form on $V\cap B_\eps$ described in Subsection~\ref{sect:Eu}.

\begin{definition}
The $G$-{\em equivariant local Euler obstruction}\index{equivariant local Euler obstruction}\index{local Euler obstruction! equivariant} of the 1-form $\omega$ on $(V,0)$ is defined by 
 \begin{equation*}
  {\rm Eu}^G(\omega;V,0)=\sum_{\overline{p}\in ({\Sing}\widetilde{\omega})/G}
  (-1)^{\dim{V}-\dim{V_{(p)}}}{\rm Eu}(V,V_{(p)})\cdot\ind(\widetilde{\omega}_{\vert V_{(p)}};V_{(p)},p)[Gp]\,,
 \end{equation*}
 where $p$ is a point of the orbit $\overline{p}=Gp$, $\ind(\cdot)$ is the usual index of a $1$-form
 on a smooth manifold.
 \end{definition}

It is not difficult to show that the equivariant local Euler obstruction is well defined (that is, that the definition does not depend on the choice of a 1-form $\widetilde{\omega}$) and its reduction under the natural reduction homomorphism $\mbox{R}^G_{\{e\}}: A(G)\to A(\{e\})=\ZZ$ gives the usual Euler obstruction of the 1-form $\omega$.

One has a global version of this notion defined either for a projective or for an affine variety, see Section~\ref{sect:Eu}. Let $V$ be a $G$-invariant affine variety in $\CC^N$ and let $\eta$ be a $G$-invariant real 1-form on $\CC^N$ which is radial at infinity (this means that it does not vanish on the vectors
$$
\sum_i \left(x_i\frac{\partial{\ }}{\partial x_i}+
y_i\frac{\partial{\ }}{\partial y_i}\right)
$$
for $\Vert z \Vert$ large enough, $z=(z_1, \ldots, z_N)$, $z_j=x_j+y_j\sqrt{-1}$) and has only isolated singular points on $V$.

\begin{definition}
The {\em equivariant global Euler obstruction}\index{equivariant global Euler obstruction}\index{global Euler obstruction! equivariant} of the affine variety $V$ is defined by
$$
{\rm Eu}^G(V):=\sum_{\overline{p}\in ({\rm Sing\,}\eta_{\rm rad})/G} {\rm I}^G_{G_p} ( {\rm Eu}^{G_p}(\eta_{\rm rad};V,p))
\in A(G).
$$
\end{definition}

The same definition makes sense for a projective (therefore compact) variety. The only difference is that one has to take an arbitrary 1-form $\eta$ with only isolated singular points.

\subsection{Real quotient singularities}\label{subsec:real_quotient}
Let a finite 
group $G$ act (linearly) on the space
$\RR^n$ (and thus on its complexification $\CC^n$). For an
analytic (real) 1-form $\omega$ on $(\RR^n, 0)$, there is defined a natural (Eisenbud--Levine--Khimshiashvili) quadratic form $B$ on 
\[
\Omega_\omega := \Omega^n_{\RR^n,0}/\omega \wedge \Omega^{n-1}_{\RR^n,0}: 
\]
see Section~\ref{subsec:algebraic_form-smooth}.
Its signature is equal to the index $\ind(\omega; \RR^n,0)$ of the 1-form $\omega$. If the 1-form $\omega$ is $G$-invariant,
its (equivariant) index $\ind^G(\omega; \RR^n,0)$ is defined as
an element of the Burnside ring $A(G)$. 
In this case
the Eisenbud--Levine--Khimshiashvili quadratic form is also
$G$-invariant and therefore its (equivariant) signature
$\sgn^G B$ is defined as an element of the ring $R_{\RR}(G)$ of real representations of the group $G$.
One can expect a relation between the equivariant signature
$\sgn^G B$ and the equivariant index $\ind^G(\omega; \RR^n,0)$
(or rather its reduction $r(\ind^G(\omega; \RR^n,0))$ under the natural homomorphism $r: A(G)\to R_{\RR}(G)$).

The most straightforward conjecture would be that the 
reduction $r(\ind^G(\omega; \RR^n,0))$ is equal to the
equivariant signature $\sgn^G B$. In~\cite{Gusein-Zade86} and also
in~\cite{Damon91}, it was explained that (for differentials of function germs) this was not the case. The reason is 
roughly speaking the following.
For a ($G$-invariant) morsification of a function germ, the usual signature
of the residue pairing can be expressed in terms of the real
critical points of the morsification, whence an equation for
the equivariant signature involves also critical points whose
complex conjugates lie in the same $G$-orbit.

A weaker conjecture can be as follows. Let $r^{(0)}: A(G) \to \ZZ$ be the group homomorphism defined by $r^{(0)}([G/H])=1$, that is, $r^{(0)}\left(\sum a_H [G/H]\right)=\sum a_H$. Let 
$B^G_\omega : \Omega^G_\omega \times
\Omega^G_\omega \to \RR$ be the restriction of the residue pairing to the
$G$-invariant part $\Omega^G_\omega$ of $\Omega_\omega$. It is a non-degenerate bilinear form as well.
It is possible to show that the image of the index
${\rm ind}^G(\omega; \CC^n,0)$ under the map $r:A(G) \to R_\CC(G)$ is equal to the class
$[\Omega^\CC_\omega]$ of the $G$-module $\Omega^\CC_\omega$: \cite{GM17}. Therefore, for the $G$-invariant
part $\left(\Omega^{\CC}_\omega\right)^G$ of $\Omega^\CC_\omega$, one has 
\begin{equation} \label{eq:r0ind}
\dim \left(\Omega^{\CC}_\omega\right)^G = r^{(0)}({\rm ind}^G(\omega; \CC^n, 0)) \, .
\end{equation}
Taking into account relations between dimensions of modules in the complex case and signatures of
quadratic forms in the real case in the Eisenbud--Levine--Khimshiashvili theory, one can conjecture that
\[
{\rm sgn}\, B_\omega^G = r^{(0)}({\rm ind}^G(\omega; \RR^n, 0)) \, .
\]
Again, in general this is not the case.

Let $W$ be the real part of the quotient $\CC^n/G$. Note that in general $W \neq \RR^n/G$. A real analytic 1-form
$\eta$ on $W$ defines a $G$-invariant analytic 1-form $\omega=\pi^\ast \eta$ on $\CC^n$ ($\pi: \CC^n \to \CC^n/G$
is the quotient map) which is real (that is real on 
$\RR^n\subset \CC^n$) and, moreover, real on $\pi^{-1}(W)$. 

One can prove the following statement (\cite{EGFAA}, \cite{EGMN}):

\begin{theorem}
For an abelian finite group $G$ and for a real analytic $G$-invariant 1-form $\omega$ one has
\begin{equation} \label{eqn:main}
{\rm sgn}\, B_\omega^G = r^{(0)}({\rm ind}^G(\omega; \pi^{-1}(W),0)) \, .
\end{equation}
\end{theorem}

\begin{remark}
 It is very probable that the statement holds for non-abelian groups as well. However, this is not proved.
\end{remark}


\printindex

\end{document}